\documentclass{amsart}
\usepackage{amssymb}
\usepackage{amsmath}
\usepackage{amsfonts}
\newtheorem{theorem}{Theorem}[section]
\newtheorem{lemma}[theorem]{Lemma}
\newtheorem{corollary}[theorem]{Corollary}

\newtheorem{prop}[theorem]{Proposition}

\newtheorem{definition}[theorem]{Definition}
\newtheorem{remark}[theorem]{Remark}

\renewcommand{\l}{\lambda}
\newcommand{\R}{\mathbb R}

\newcommand{\p}{\partial}

\newcommand{\ta}{{\tilde{a}}}

\newcommand{\tc}{{\tilde{c}}}
\newcommand{\tS}{{\tilde S}}

\newcommand{\tw}{{\tilde{w}}}
\newcommand{\tA}{{\tilde{A}}}
\newcommand{\tX}{{\tilde{X}}}
\newcommand{\tphi}{{\tilde{\phi}}}
\begin{document}
\title{Gradient NLW on curved background in $4+1$ dimensions}

\author{Dan-Andrei Geba and Daniel Tataru}

\address{Department of Mathematics, Hylan Building \\University of Rochester, Rochester, NY 14627}
\address{Department of Mathematics, Evans Hall \\University of California
at Berkeley, Berkeley, CA 94720-3840}
\date{}

\begin{abstract}
  We obtain a sharp local well-posedness result for the Gradient
  Nonlinear Wave Equation on a nonsmooth curved background.  In the
  process we introduce variable coefficient versions of Bourgain's
  $X^{s,b}$ spaces, and use a trilinear multiscale wave packet
  decomposition in order to prove a key trilinear estimate.
\end{abstract}

\maketitle
\section{Introduction}

In this article we are investigating the issue of local well-posedness
for a variable coefficient semilinear wave equation in $4+1$ dimensions.
To describe the context and motivate the interest in  our problem
we introduce three related equations. We begin with a generic
gradient NLW equation in $\R^{n+1}$,
\begin{equation}
\Box u\,=\,\Gamma(u)(\nabla u)^2
\label{sem}
\end{equation}
with the nonlinearity
\[
\Gamma(u) (\nabla u)^2\,=\,q^{ij}(u) \p_i u\, \p_j u
\]
where $q^{ij}$ are smooth functions and the standard summation
convention is used.

Then we move on to a similar equation but on a curved background,
\begin{equation}
\Box_g u\,=\,\Gamma(u)(\nabla u)^2
\label{eq}
\end{equation}
with $\Box_{g}\,=\,g^{ij}\,\partial_i\partial_j$, where the
summation occurs from $0$ to $n$ and the index $0$ stands for the
time variable. To insure hyperbolicity we assume that the matrix
$g^{ij}$ has signature $(1,n)$ and the time level sets $x_0=const$
are space-like, i.e. $g^{00}>0$. In effect to simplify some of the
computations we make the harmless assumption $g^{00}=1$.

Finally, we consider a corresponding quasilinear equation
\begin{equation}
\Box_{g(u)} u\,=\,\Gamma(u)(\nabla u)^2
\label{qua}
\end{equation}
with similar assumptions on the matrix $g$.

In all three cases we are interested in the local well-posedness of
the Cauchy problem in Sobolev spaces $H^s(\R^n)\times H^{s-1}(\R^n)$
with initial data
\begin{equation}
u(0,x)\,=\,u_0(x),\qquad \partial_t u(0,x)\,=\,u_1(x)
\label{id}
\end{equation}
The first equation \eqref{sem} is the best understood so far, and
is known to be locally well-posed for $s$ in the range
\[
s>\max\{\frac{n}{2},\frac{n+5}{4}\}
\]
This range is sharp. The $\frac{n}{2}$ obstruction comes from
scaling, while the $\frac{n+5}{4}$ is related to concentration
along light rays, see Lindblad \cite{MR1375301}. The proof of the
positive result is fairly straightforward in dimension $2+1$ and
$3+1$, where it suffices to rely on the Strichartz estimates. In
$4+1$ dimensions this no longer works and one needs to use instead
the $X^{s,\theta}$ spaces, see Foschi-Klainerman \cite{MR1755116}.
These are multiplier weighted $L^2$ spaces associated to the wave
operator as the Sobolev spaces $H^s$ are connected to the Laplace
operator $\Delta$, see Klainerman-Machedon \cite{MR1420552}:
\begin{equation}
\| u \|_{X^{s,\theta}} = \| (1+|\xi|^2)^{\frac s2}\cdot
(1+||\tau|-|\xi||^2)^{\frac \theta2}\cdot |\hat{u}(\tau,\xi)|
\|_{L^2}
\label{standardxs}\end{equation}
where $\hat u=\hat{u}(\tau,\xi)$ is the space-time Fourier
transform of function $u=u(t,x)$.  Finally, in the most difficult
case, $n\geq 5$, this was proved by Tataru \cite{MR1739207}, using
a suitable modification of the $X^{s,\theta}$ spaces, needed in
order to control the interaction of high and low frequencies in
the multiplicative estimates.

For the quasilinear problem (\ref{qua}) the sharp result is only known
to hold in dimensions $n=2,3$. This was proved by Smith-Tataru
\cite{MR2178963} (see also Lindblad's counterexample
\cite{MR1666844}). The argument there still requires the use of
Strichartz estimates. These are derived from a wave packet parametrix
construction for a wave equation with very rough coefficients, which
in turn is obtained via a very delicate analysis of the Hamilton flow.
A different proof of this result in the special case of the Einstein
vacuum equation was independently obtained by Klainerman-Rodnianski
\cite{MR1885093}, \cite{MR2180401}, \cite{MR2052472}.  In dimensions
$n\geq 4$ it is still unclear which is the optimal threshold, the best
results so far being contained in the above mentioned paper of
Smith-Tataru~\cite{MR2178963} and in an earlier one,
Tataru~\cite{MR1887639}:
\[
\aligned &n=4,5\quad \quad s>\frac{n}{2}+\frac{1}{2}\\
&n\geq 6\qquad \quad s>\frac{n}{2}+\frac{2}{3}
\endaligned
\]
In the same direction but somewhat closer in spirit to the present
paper is Bahouri and Chemin's work \cite{MR2082388,MR2003418}.
The equation considered there is still quasilinear, but the main
estimates are frequency localized versions of the Strichartz
estimates for the wave equation on a rough background.

As an intermediate step toward understanding the higher dimensional
quasilinear problem, we consider here the semilinear problem on a
curved background and we prove the sharp result:

\begin{theorem}
Let $n=4$ and assume that the coefficients $g^{ij}$ satisfy
$\partial^2 g \in L^2 L^\infty$.  Then the Cauchy problem
\eqref{eq}, \eqref{id} is locally well-posed in $H^s\times
H^{s-1}$ for $s>\frac{9}{4}$.
\label{main}
\end{theorem}

Here well-posedness is understood in the strongest sense, i.e.  the
solutions have Lipschitz dependence on the initial data and they exist
on a time interval which only depends on the size of the initial data.

One contribution of the present paper is to introduce variable
coefficient versions of the $X^{s,b}$ spaces, study their properties
and obtain the corresponding Strichartz type embeddings.
However, the main novelty, contained in the last two sections,
is a new method, based on a trilinear wave packet decomposition,
to prove a key trilinear bound which cannot be obtained directly
from the Strichartz estimates.

The first step in the proof is to reduce the problem to the case
when the initial data is small, using scaling and the finite speed
of propagation. This is a routine argument for which we refer the
reader to \cite{MR2178963}.  Once we know that the initial data is
small, we can fix the time interval and set it to $[-1,1]$. This
will be the case throughout the rest of the paper.

To solve the problem for small data we use a fixed point argument.
Let  $S(u_0,u_1)$ and $\Box_{g}^{-1}$ be respectively the
homogeneous and inhomogeneous solution operators
\begin{align}
&\Box_{g}S(u_0,u_1)\,=\,0,\quad S(u_0,u_1)(0)\,=\,u_0,\quad
\partial_t S(u_0,u_1)(0)\,=\,u_1
\label{hom}\\
&\Box_{g}(\Box_{g}^{-1}H)\,=\,H,\qquad
(\Box_{g}^{-1}H)(0)\,=\,0,\qquad
\partial_t(\Box_{g}^{-1}H)(0)\,=\,0
\label{inhom}
\end{align}
Then a solution $u$ for (\ref{eq}) in $[-1,1]$ is
also a fixed point for the functional
\begin{equation}
F(u)\,= S(u_0,u_1) + \Box_{g}^{-1}( \Gamma(u)(\nabla
u)^2)
\label{funct}
\end{equation}

In order to apply a fixed point argument for $F$ we need to find
two Banach spaces $X$ and $Y$ for which the following mapping
properties hold:
\begin{align}
&\,\| S(u_0,u_1) \|_X \lesssim \|(u_0,u_1)\|_{H^s\times H^{s-1}}
\label{h}\\
&\,\|\Box_g^{-1}H\|_X \lesssim \|H\|_Y
\label{i}\\
&\,\|u \cdot w\|_X \lesssim \|u\|_X \|w\|_X
\label{xx}\\
&\|\Gamma(u)\|_{X} \lesssim C(\|u\|_{L^\infty}) (1 + \|u\|_{X}^5)
\label{moser} \\
&\,\|u \cdot w\|_Y \lesssim \|u\|_X \|w\|_Y
\label{xy} \\
&\,\|\nabla v\cdot \nabla w\|_Y\lesssim  \|v\|_X\cdot \|w\|_X
\label{n}
\end{align}
where $C=C(\|u\|_{L^\infty})$ is a constant that depends solely on
$\|u\|_{L^\infty}$. In the flat case (\ref{sem}), for dimension
$n=4$, one can make this argument work by choosing
\[
X\,=\,X^{s,\theta}\qquad Y\,=\,X^{s-1,\theta-1}
\]
with
\begin{equation}
s\,=\,\theta\,+\frac 32\qquad \theta\,>\,\frac 34 \label{st}
\end{equation}

For our problem the challenge is twofold: first we need to find
suitable variable coefficient versions for the $X^{s,\theta}$
spaces and then, in this new context, prove the corresponding
estimates (\ref{h})-(\ref{n}).

Such spaces were previously introduced by Tataru \cite{MR1391526},
where they are used in the context of a unique continuation
problem. There, for a hyperbolic operator $P$  one  defines
\[
X^{s,0}\,=\,H^s,\quad X^{s,1}\,=\,\{u\in H^s| Pu\in H^{s-1}\}
\]
Then all the other spaces are defined through interpolation and
duality.

In this article we choose to follow a different path based on dyadic
decompositions with respect to the spatial frequency and the distance
to the characteristic cone. Likely one should be able to prove that
the two approaches are equivalent, but we choose not to pursue this
here.

Our article is structured as follows. In the next section we define
the $X^{s,\theta}$ spaces and prove that they satisfy the linear
estimates \eqref{h}, \eqref{i}. Our definition of the $X^{s,\theta}$
is slightly different from the standard one \eqref{standardxs} in the constant
coefficient case. Precisely, in the constant coefficient case
our definition gives
\begin{equation}
\| u \|_{X^{s,\theta}} \approx  \| (1+|\xi|^2)^{\frac s2}
(1+||\tau|-|\xi||^2)^{\frac \theta2}\cdot \hat{u}(\tau,\xi)
\|_{L^2} + \|\Box u\|_{L^2_t H^{s+\theta-2}_x}, \quad 0 < \theta <
1 \label{newxs}\end{equation} and one can see that the second term
above alters the behavior at high modulations $|\tau| \gg |\xi|$.
Correspondingly, for negative $\theta$ we
 have
\begin{equation}
\| u \|_{X^{s,\theta}} \approx  \| (1+|\xi|^2)^{\frac s2}
(1+||\tau|-|\xi||^2)^{\frac \theta2}\cdot \hat{u}(\tau,\xi)
\|_{L^2} + \| u\|_{L^2_t H^{s+\theta}_x}, \quad -1 < \theta < 0
\label{newxs1}\end{equation}
This change is consistent with scaling and simplifies somewhat
the study of high modulation interactions.

In Section~\ref{secse} we discuss the
Strichartz estimates for $\Box_g$, which translate into embeddings
for the $X^{s,\theta}$ spaces. These turn out to suffice for the
proof of the algebra properties \eqref{xx}-\eqref{xy} and for the
high-high frequency interactions in \eqref{n}.

The difficult part is  to study the high-low frequency interactions
in \eqref{n}. For this we first take advantage of  the duality relation
\begin{equation}
(X^{s,\theta}+ L^2 H^{s+\theta})'\,=\,X^{-s,-\theta} \label{dual}
\qquad s \in \R, \quad 0 < \theta < \frac 12
\end{equation}
This is consistent with \eqref{newxs} and \eqref{newxs1}.  Using
this duality, after factoring out high modulation interactions,
the bound (\ref{n}) is transformed into the trilinear estimate:
\begin{equation}
\left|\int u\cdot v\cdot w\,\,dx \,dt\right|\lesssim \|u\|_{X^{1-s,1-\theta}}
\|v\|_{X^{s-1,\theta}} \|w\|_{X^{s-1,\theta}}
\label{trilin}
\end{equation}
with $(s,\theta)$ verifying (\ref{st}). The last section of the paper
is devoted to proving this bound. The argument is based on a
multiscale trilinear  wave packet decomposition for linear waves.

\section{The $X^{s,\theta}$ spaces}

We first introduce  Littlewood-Paley decompositions. As a general
rule, all frequency localizations in the sequel are only with respect
to the spatial variables. There is a single exception to this. Precisely,
the coefficients $g^{ij}$ are truncated using space-time
multipliers. In order for these truncations to work, we need for these
coefficients to be defined globally in time. Hence we assume they
have been extended to functions with similar properties in all of
$\R^{n+1}$.

Let $\phi$ be a smooth function supported in $\{\frac12 \leq |\xi|
\leq 2 \}$ with the property that
\[
1 = \sum_{j=-\infty}^\infty \phi(2^{-j} \xi)
\]
We consider a spatial Littlewood-Paley decomposition,
\[
1 = \sum_{\l=1}^\infty S_{\l}(D_x)
\]
where for dyadic $\l > 1$ we have
\[
\qquad S_{\l}(\xi) = \phi\left(\frac{\xi}{\l}\right)
\]
while $S_1$ incorporates the low frequency contribution in
$\{|\xi|\leq 1\}$. Set
\[
S_{<\l} =  \sum_{\mu=1}^{\frac\l2} S_\mu
\]
We will also use spatial multipliers $\tS_\l$ with slightly larger
support, with $S_\l \tS_\l = S_\l$.  We say that a function $u$ is
localized at frequency $\l$ if its Fourier transform is supported in
the annulus $\{ \frac\l 8 \leq |\xi| \leq 8\l\}$.

For the paradifferential type calculus we
also need to truncate the coefficients of $\Box_g$ in frequency.
Given $\Box_g$ in (\ref{eq}) we define the modified operators
\[
\Box_{g_{<\l}} =(S_{<\l}(D_x,D_t) g^{\alpha\beta}) \partial_\alpha
\partial_\beta
\]
In the sequel we omit the space and time variables in our function
space notations, i.e. $L^p:= L^p_{x,t}$, $L^2 H^s:= L^2_t H^s_x$, $L^p
L^q:= L^p_t L^q_x$, etc. We are ready now to define our spaces:

\begin{definition}
Let $\theta\in (0,1)$ and $s \in \R$. Then $X^{s,\theta}$
is the space of functions $u \in L^2(-1,1; H^s(\R^n))$
for which the following norm is finite:
\begin{equation}
\|u\|^2_{X^{s,\theta}}\,=\,\inf \left \{ \sum_{\l=1}^\infty
\sum_{d=1}^\lambda \|u_{\lambda,d}\|^2_{X_{\l,d}^{s,\theta}} ;\ u
=\sum_{\l=1}^\infty\sum_{d=1}^{\l}S_{\l}
u_{\l,d}\right\}
\label{tp}
\end{equation}
where $\l$, $d$ take dyadic values and
\begin{equation}
\|u_{\l,d}\|_{X_{\l,d}^{s,\theta}}^2 =
\l^{2s}d^{2\theta}\|u_{\l,d}\|^2_{L^2} +
\l^{2s-2}d^{2\theta-2}\|\Box_{g_{<{\sqrt\l}}}u_{\l,d}\|^2_{L^2}
\label{cl}
\end{equation}

We also define the space $X^{s-1,\theta-1}$ of functions for which
the following norm is finite:
\begin{equation}
\begin{split}
\| f\|_{X^{s-1,\theta-1}}^2 = \inf \left\{ \|f_0\|_{L^2 H^{s-1}}^2
+ \sum_{\l=1}^\infty\right.& \sum_{d=1}^\lambda
\|f_{\lambda,d}\|^2_{X_{\l,d}^{s,\theta}} ;\ \\ & \left. f =f_0+
\sum_{\l=1}^\infty \sum_{d=1}^{\l}\Box_{g_{<{\sqrt\l}}} S_{\l}
f_{\l,d} \right\}
\end{split}
\label{tn}
\end{equation}
\label{xst}
\end{definition}
\begin{remark}
  Intuitively $d$ stands for the modulation of the
  $u_{\l,d}$ piece. Indeed, in the constant coefficient case
one can easily see that $u_{\lambda,d}$  mainly contributes
to $u$ in the region where $||\tau|-|\xi|| \approx d$.
The condition $1 \leq d$ is related to the spatial localization on the
unit scale in our problem. The condition $d \leq \lambda$
reflects the fact that at high modulation we use a simpler structure,
see e.g. \eqref{newxs}, \eqref{newxs1}.
\end{remark}
\begin{remark}
  The cutoff at frequency less than $\sqrt\l$ for the coefficients
  $\Box_g$ is related to the regularity of the coefficients,
  $\partial^2 g \in L^2 L^\infty$.  This implies that
  $\Box_{g_{\geq{\sqrt\l}}}u_{\l,d}$ is an allowable error term.\\
\end{remark}

We begin our analysis of  the $X^{s,\theta}$ spaces with a simple observation,
namely that without any restriction in generality one can assume that
the functions $u_{\l,d}$ and $f_{\l,d}$ in Definition~\ref{xst} are
localized at frequency $\l$. Precisely, we have the stronger result:

\begin{lemma}
The following estimate holds:
\begin{equation}
\l^{s-1}d^{\theta}\|\nabla S_\l v\|_{L^2} +
\l^{s-1}d^{\theta-1}\|\Box_{g_{<{\sqrt\l}}}S_\l v \|_{L^2}
\lesssim \|v\|_{X_{\l,d}^{s,\theta}} \label{strongst}
\end{equation}
\end{lemma}
\begin{proof}
We first bound the time derivatives of $v$ in negative Sobolev spaces,
\begin{equation}
\l^{s} d^{\theta} (\| \p_t^2 v\|_{L^2 (H^{-2}+\l^2 L^2)} + \|\p_t
v\|_{L^2(H^{-1}+\l L^2)}  ) \lesssim \|v\|_{X_{\l,d}^{s,\theta}}
\label{lowv}
\end{equation}
This follows by Cauchy-Schwartz from the interpolation inequality
\[
\|\p_t v\|_{L^2 ( H^{-1}+\l L^2)}^2 \lesssim  (\| \p_t^2 v\|_{L^2 (H^{-2}+\l^2L^2)}
+\|v\|_{L^2})\|v\|_{L^2}
\]
combined with the bound
\[
\| \p_t^2 v\|_{L^2 (H^{-2}+\l^2 L^2)} \lesssim \l^{-2} \|
\Box_{g_{<{\sqrt\l}}} v\|_{L^2} + \|\p_t v\|_{L^2 ( H^{-1}+\l
L^2)} + \|v\|_{L^2}
\]
To prove this last estimate we only use the $L^\infty$ regularity of  $g$
together with the condition $g^{00} = 1$. Then we need  the
fixed time bounds
\[
\| g_{<{\sqrt{\l}}} \p_x \p_t v\|_{H^{-2}+\l^2 L^2} \lesssim
\|\p_t v\|_{H^{-1} + \l L^2}
\]
\[
\| g_{<{\sqrt{\l}}} \p_x^2 v\|_{ H^{-2}+\l^2 L^2} \lesssim
\|v\|_{L^2}
\]
They are similar, so we only discuss the second one. We write
\[
 g_{<{\sqrt{\l}}} \p_x^2 v = \p_x^2  (g_{<{\sqrt{\l}}} v) - 2\p_x
 (\p_x g_{<{\sqrt{\l}}} v) +  \p_x^2 g_{<{\sqrt{\l}}} v
\]
and use the uniform bounds
\[
|g_{<{\sqrt{\l}}}| \lesssim 1, \qquad |\p_x g_{<{\sqrt{\l}}}|
\lesssim \l,
 \qquad |\p_x^2 g_{<{\sqrt{\l}}}| \lesssim \l^2
\]
This concludes the proof of \eqref{lowv}.

The first term in \eqref{strongst} is directly bounded using
\eqref{lowv}. For the second it suffices to prove  the commutator estimate
\begin{equation}
\|[\Box_{g_{<{\sqrt\l}}},S_\l] v\|_{L^2} \lesssim \l \|v\|_{L^2}+
\|\p_t v\|_{ L^2} \label{smallgcom}
\end{equation}
We have
\[
[\Box_{g_{<{\sqrt\l}}},S_\l] = [g_{<{\sqrt{\l}}},S_\l] \p_t \p_x +
[g_{<{\sqrt{\l}}},S_\l] \p_x^2
\]
and the commutators are localized at frequency $\l$ so the spatial
derivatives only contribute factors of $\l$. Hence \eqref{smallgcom}
follows from the standard commutator estimate
\[
\|[g_{<{\sqrt{\l}}},S_\l]\|_{L^2 \to L^2} \lesssim \l^{-1}
\|\nabla g\|_{L^\infty}
\]
\end{proof}

Applying the above Lemma with $S_\l$ replaced by $\tS_\l$ we obtain

\begin{corollary}
One can replace the $X^{s,\theta}_{\l,d}$ norm in the definition
of $X^{s,\theta}$ and $X^{s-1,\theta-1}$ by the norm
\[
\| v\|_{\tilde{X}^{s,\theta}_{\l,d}} = \l^{s-1}d^{\theta}\|\nabla
v\|_{L^2} + \l^{s-1}d^{\theta-1}\|\Box_{g_{<{\sqrt\l}}} v\|_{L^2}
\]
\label{tx}
\end{corollary}

For the proof of the duality relation \eqref{dual} it is convenient
to work with a selfadjoint operator. Thus we consider the selfadjoint
counterpart $\tilde{\Box}_g$ of $\Box_g$
\[
\tilde{\Box}_g = \partial_i g^{ij} \partial_j
\]
Then for $v$ localized at frequency $\lambda$ we commute and estimate
the frequency localized difference
\[
\| \tilde \Box_{g_{<{\sqrt\l}}} v -  \Box_{g_{<{\sqrt\l}}} v\|_{L^2}
\lesssim \|\nabla v\|_{L^2}
\]
This leads directly to

\begin{corollary}
One can replace the $\Box_{g_{<{\sqrt\l}}} $ operator in the
definition of $X^{s,\theta}$ and $X^{s-1,\theta-1}$ by the similar
operator in divergence form $\tilde\Box_{g_{<{\sqrt\l}}} $. \label{self}
\end{corollary}

As a consequence of the second part of \eqref{strongst} we
have

\begin{corollary}
The following embedding holds for $-1 < \theta < 0$:
\[
X^{s,\theta} \subset L^2 H^{s+\theta}
\]
\label{sw}
\end{corollary}
Another use of this is to establish energy estimates. A direct
application of energy estimates for the wave equation yields
the bound
\[
\| \nabla v\|_{L^\infty L^2}^2 \lesssim \|\nabla v\|_{L^2}^2 +
\|\nabla v\|_{L^2} \| \Box_{g} v\|_{L^2}
\]
This leads to
\begin{equation}
\l^{s-1} d^{\theta-\frac12} \| \nabla S_\l v\|_{L^\infty L^2}
\lesssim \| v\|_{\tilde{X}^{s,\theta}_{\l,d}} \label{litwo}
\end{equation}
Going back to Definition~\ref{xst}, this implies

\begin{corollary}
Assume that $\theta > \frac12$. Then
\begin{equation}
\|u\|_{L^\infty H^s} + \|u_t\|_{L^\infty H^{s-1}} \lesssim \|u \|_{X^{s,\theta}}
\end{equation}
\label{trace}
\end{corollary}

To prove the estimates \eqref{h} and \eqref{i} in the context of the
$X^{s,\theta}$ spaces we  need to switch from the frequency truncated
coefficients to the full coefficients $g^{ij}$.  The tool needed to do
that is contained in the following:

\begin{lemma}
Assume that $0 \leq s \leq 3$. Then the following fixed time
estimate holds:
\begin{equation}
\sum_{\l=1}^\infty \lambda^{2(s-1)}\|\tS_\l ({g_{> \sqrt{\l}}} u)
\|_{L^2}^2 \lesssim (M (\|\p^2 g\|_{L^\infty}))^2 \| u\|_{
H^{s-2}}^2 \label{largeg}
\end{equation}
where $M$ stands for the maximal function with respect to time.
We also have the dual estimate
\begin{equation}
\| \sum_{\l=1}^\infty  {g_{> \sqrt{\l}}} \tS_\l f_\l
\|_{H^{2-s}}^2 \lesssim \sum_{\l=1}^\infty \lambda^{2(1-s)}\| f_\l
\|_{L^2}^2 \label{largegdual}
\end{equation}
\end{lemma}

\begin{proof}
We take a Littlewood-Paley decomposition of both factors,
\[
\tS_\l  (g_{> \sqrt{\l}} u)  = \sum_{\mu=1}^\infty \sum_{\nu =
\sqrt{\l}}^\infty \tS_\l( g_\nu  u_\mu)
\]
The $(\mu,\nu)$ term is nonzero only in the following situations:

(i) $\nu \ll \l$, $\mu \approx \l$. Then we estimate
\[
\| \tS_\l (g_\nu  u_\mu)\|_{L^2} \lesssim \|g_\nu\|_{L^\infty}
\|u_\mu\|_{ L^2} \lesssim  \nu^{-2}  M (\|\p^2 g\|_{L^\infty})
\|u_\mu\|_{ L^2}
\]
and use the square summability with respect to $\l$ together with the relation
$\nu^{-2} \lesssim  \l^{-1}$.

(ii) $\nu \approx \l$, $\mu \ll \l$. Then
\[
\| \tS_\l (g_\nu  u_\mu)\|_{L^2} \lesssim \|g_\nu\|_{L^\infty}
\|u_\mu\|_{ L^2} \lesssim \l^{-2} M (\|\p^2 g\|_{L^\infty})
\|u_\mu\|_{ L^2}
\]
This is tight only when $s=3$ and $\mu=1$, otherwise there is a gain
which insures the summability in $\l$, $\mu$.

(iii) $\nu \approx \mu \gtrsim \l$. Then
\[
\| \tS_\l (g_\nu  u_\mu)\|_{L^2} \lesssim \|g_\nu\|_{L^\infty}
\|u_\mu\|_{ L^2} \lesssim  \mu^{-2} M (\|\p^2 g\|_{L^\infty}) \|
u_\mu\|_{ L^2}
\]
This is always stronger than we need. The proof of the lemma is
concluded.
\end{proof}

We now establish some simple properties of  the linear equation
\begin{equation}
\Box_g u = f, \qquad u(0) = u_0, \qquad u_t(0) = u_1.
\label{leq}
\end{equation}
Then

\begin{lemma}
The linear equation \eqref{leq} is well-posed in $H^s \times
H^{s-1}$ for $0 \leq  s \leq 3$.
\label{lwp}
\end{lemma}
The proof follows easily from energy estimates, see
\cite{MR2153517}.

We use this to prove \eqref{h}, namely

\begin{lemma}
Assume that $0 \leq  s \leq 3$ and $\theta > 0$. Then the solution $u$
to \eqref{leq} verifies
\[
\| u \|_{X^{s,\theta}} \lesssim \|u_0\|_{H^s} +
\|u_1\|_{H^{s-1}} + \|f\|_{L^2 H^{s-1}}
\]
\label{nhom}
\end{lemma}

\begin{proof}
We decompose the solution $u$ as
\[
u = \sum_{\l=1}^\infty S_\lambda  \tS_\l u
\]
and think of this as a part of the sum in \eqref{tp} which corresponds
to $d=1$. Then
\[
\begin{split}
\|u\|_{X^{s,\theta}}^2 &\lesssim \sum_{\l=1}^\infty \|\tS_\l
u\|_{X^{s,\theta}_{\lambda,1}}^2 \\
& \approx  \sum_{\l=1}^\infty \l^{2s} \| \tS_\l u\|^2_{L^2} +
\lambda^{2(s-1)} \|\Box_{g_{<{\sqrt\l}}}  \tS_\l u\|^2_{L^2}
\\ &\lesssim \|  u\|^2_{L^2 H^s} +  \sum_{\l=1}^\infty
\lambda^{2(s-1)} \|\Box_{g_{<{\sqrt\l}}}  \tS_\l u - \tS_\l \Box_g
u \|^2_{L^2} + \|f\|^2_{L^2 H^{s-1}}
\end{split}
\]
The first term is easily controlled by energy estimates. The
second is decomposed as follows:
\[
\Box_{g_{<{\sqrt\l}}}  \tS_\l u - \tS_\l \Box_g u =
[\Box_{g_{<{\sqrt\l}}},\tS_\l] u - \tS_\l \Box_{g_{>{\sqrt\l}}} u
\]
For the commutator we use the fixed time bound \eqref{smallgcom}
along with square summability in $\lambda$. The second part is
controlled by \eqref{largeg}.

\end{proof}

The result in the next Lemma implies the estimate \eqref{i} for the
spaces $X,Y$:

\begin{lemma}
Assume that $0 \leq s \leq 3$ and $\frac12 < \theta < 1$.
Then the operator $\Box_g^{-1}$ has the mapping property
\[
\Box_g^{-1}: X^{s-1,\theta-1} \to X^{s,\theta}
\]
\end{lemma}

\begin{proof}
Let $f \in X^{s-1,\theta-1}$. We use the representation in
\eqref{tn},
\[
f = f_0 + \sum_{\l=1}^\infty \sum_{d=1}^{\l}\Box_{g_{<{\sqrt\l}}}
S_{\l} f_{\l,d}
\]
By Definition~\eqref{xst} the function
\[
u = \sum_{\l=1}^\infty \sum_{d=1}^{\l} S_{\l} f_{\l,d}
\]
belongs to $X^{s,\theta}$. The difference $v = u - \Box_g^{-1} f$
solves
\[
\Box_g v = \Box_g u - f, \qquad v(0) = u(0), \qquad v_t(0) = u_t(0)
\]
To estimate it we use Lemma~\ref{nhom}. The initial data is
controlled due to Corollary~\ref{trace}, so it remains to bound
the inhomogeneous term in $L^2 H^{s-1}$. Thus we need to show that
\begin{equation}
\| \sum_{\l=1}^\infty  \sum_{d=1}^{\l} \Box_{g_{>{\sqrt\l}}}
S_{\l} f_{\l,d} \|_{L^2H^{s-1}}^2 \lesssim \sum_{\l,d}
\|f_{\l,d}\|_{X^{s,\theta}_{\l,d}}^2
\end{equation}
Considering the trace regularity result in Corollary~\ref{trace}
this would follow from
\[
\| \sum_{\l=1}^\infty \Box_{g_{>{\sqrt\l}}} S_{\l} f_{\l}
\|_{L^2H^{s-1}}^2 \lesssim \sum_{\l}  \|\nabla f_{\l}\|_{L^\infty
H^{s-1}}^2, \qquad f_\l = \sum_{d=1}^\l f_{\l,d}
\]
which in turn is a consequence of the fixed time bound \eqref{largegdual}.

\end{proof}

We finish this section by proving a key duality relation between
$X^{s,\theta}$ spaces with positive, respectively negative $\theta$.

\begin{lemma}
For $ 0 < \theta < \frac12$ we have the duality relation
\begin{equation}
X^{-s,-\theta}= (X^{s,\theta} + L^2 H^{s+\theta})'
\label{dual1} \end{equation}
\end{lemma}
\begin{proof}
a) We first show that
\[
X^{-s,-\theta} \subset (X^{s,\theta} + L^2 H^{s+\theta})'
\]
From Corollary~\ref{sw} we obtain $
X^{-s,-\theta} \subset (L^2 H^{s+\theta})' $.
It remains to prove the bound
\[
\left|\int u\cdot f\,\, dx\,dt\right|\lesssim
\|u\|_{X^{s,\theta}}\,\|f\|_{X^{-s,-\theta}}
\]
We consider Littlewood-Paley decompositions of $u$ and $v$
as in Definition \ref{xst},
\[
u =  \sum_{\l=1}^\infty  \sum_{d=1}^{\l}  S_\lambda u_{\lambda,d},
 \qquad f = f_0 + \sum_{\l=1}^\infty  \sum_{d=1}^{\l}  \tilde \Box_{g_{<\sqrt{\lambda}}}
S_\lambda f_{\lambda,d}
\]
with $\tilde \Box_{g_{<\sqrt{\lambda}}}$ in divergence form, see
Corollary~\ref{self}.  The summation with respect to $\l$ is
essentially diagonal therefore it follows by orthogonality. To handle
the $d$ summation it suffices to obtain the off-diagonal decay
\[
\begin{split}
\left | \int S_\l u_{\lambda,d_1} \cdot\tilde\Box_{g_{<\sqrt{\lambda}}}
S_\l
  f_{\lambda,d_2} dx dt\right|  \lesssim &
\min\left\{\left(\frac{d_2}{d_1}\right)^\theta,
 \left (\frac{d_1}{d_2}\right)^{\frac12-\theta}\right\}
\\ & \|u_{\lambda,d_1}\|_{X^{s,\theta}_{\lambda,d_1}
}\,\|f_{\lambda,d_2}\|_{X^{1-s,1-\theta}_{\lambda,d_2}}
\end{split}
\]
If $d_2 < d_1$ then this follows directly from \eqref{cl} and
\eqref{strongst}.  Otherwise we integrate by parts
\[
\begin{split}
 \int S_\l u_{\lambda,d_1}\cdot \, &\tilde\Box_{g_{<\sqrt{\lambda}}} S_\l f_{\lambda,d_2} dx
dt = \int\tilde \Box_{g_{<\sqrt{\lambda}}}  S_\l u_{\lambda,d_1} \cdot
S_\l f_{\lambda,d_2} dx dt \\  + & \left. \int  (S_\l
u_{\lambda,d_1}\cdot g^{0\alpha}_{<\sqrt{\l}}\partial_\alpha  S_\l
f_{\lambda,d_2} -
 g^{0\alpha}_{<\sqrt{\l}}\partial_\alpha  S_\l u_{\lambda,d_1} \cdot S_\l f_{\lambda,d_2}) dx \right|_{-1}^1
\end{split}
\]
For the first term we use \eqref{strongst} and \eqref{cl}. For the
second we use the trace regularity result in \eqref{litwo}.

b) We now show that
\[
(X^{s,\theta}+ L^2 H^{s+\theta})' \subset X^{-s,-\theta}
\]
Let $T$ be a bounded linear functional on $X^{s,\theta} + L^2
H^{s+\theta}$.  Due to the second term we can identify $T$ with a
function $u \in L^2 H^{-s-\theta}$.

On the other hand, we can apply it to functions $v \in X^{s,\theta}$
of the form
\[
v =  \sum_{\l=1}^\infty  \sum_{d=1}^{\l} S_\lambda v_{\lambda,d}
\]
Then we must have the bound
\[
|Tv|^2  \lesssim \|v\|_{X^{s,\theta}}^2  \lesssim \sum_{\l,d}
\|v_{\lambda,d}\|^2_{X^{s,\theta}_{\lambda,d}} \lesssim
\sum_{\l,d} \left( \l^{2s}d^{2\theta}\|v_{\l,d}\|^2_{L^2} +
\l^{2s-2}d^{2\theta-2}\|\tilde
\Box_{g_{<{\sqrt\l}}}v_{\l,d}\|^2_{L^2}\right)
\]
Given the definition of the $X^{s,\theta}_{\lambda,d}$ norms, using
succesively the Hahn-Banach theorem and Riesz's theorem it follows
that we can find functions $f_{\lambda,d}$ and $h_{\lambda,d}$ with
\begin{equation}
 \sum_{\l=1}^\infty  \sum_{d=1}^{\l} \lambda^{-2s} d^{-2\theta}
 \|f_{\lambda,d}\|_{L^2}^2
+   \lambda^{2(1-s)} d^{2(1-\theta)} \|h_{\lambda,d}\|_{L^2}^2 = M < \infty
\label{M}\end{equation}
so that
\[
Tv = \sum_{\l=1}^\infty \sum_{d=1}^{\l} \int f_{\lambda,d}\,
v_{\lambda,d} + h_{\lambda,d}\cdot \tilde \Box_{g_{<
\sqrt\lambda}}v_{\lambda,d}\, dx dt
\]
In particular this must hold for $v$ of the form $v = S_\lambda
v_{\lambda,d}$,
\[
\int u \, S_\lambda  v_{\lambda,d} dx dt  = \int f_{\lambda,d}\,
v_{\lambda,d} + h_{\lambda,d}\cdot \tilde \Box_{g_{<
\sqrt\lambda}}v_{\lambda,d} \, dx dt
\]
For each $\lambda,d$ this yields
\[
S_\lambda u = f_{\lambda,d} + \tilde \Box_{g_{< \sqrt\lambda}} h_{\lambda,d}
\]
Then we can represent $S_\lambda u$ in the form
\begin{equation}
S_\lambda u = f_{\lambda,1} + \sum_{d=1}^{\frac\lambda 2}
\tilde \Box_{g_{< \sqrt\lambda}} u_{\lambda,d} + \tilde \Box_{g_{<
\sqrt\lambda}} h_{\lambda,\lambda} \qquad u_{\lambda,d} =
h_{\lambda,d} -h_{\lambda,2d}
\label{py}
\end{equation}
This yields for $u$ the representation
\begin{equation}
u = \sum_{\lambda = 1}^\infty \tilde S_\lambda \left(  f_{\lambda,1} + \sum_{d=1}^{\frac\lambda 2}
\tilde \Box_{g_{< \sqrt\lambda}} u_{\lambda,d} + \tilde \Box_{g_{<
\sqrt\lambda}} h_{\lambda,\lambda}\right)
\label{pya}\end{equation}
This is very close to but not exactly the form in \eqref{tn}. However
the multipliers $\tilde S_\lambda$ can be easily replaced  by
$S_\lambda$ by reapplying the Paley-Littlewood decomposition on the
right, and then $S_\lambda$ can be commuted to the right of $\tilde \Box_{g_{<
\sqrt\lambda}}$ due to the Corollary~\ref{tx} and the commutator
bound \eqref{smallgcom}. Hence we have
\[
\| u\|_{X^{-s,-\theta} }^2 \lesssim \sum_{\lambda = 1}^\infty
\left(\lambda^{-2s} \|f_{\lambda,1} \|_{L^2}^2 +
\sum_{d=1}^{\lambda/2}
 \|  u_{\lambda,d}\|_{X_{\lambda,d}^{1-s,1-\theta}}^2 +
 \| h_{\lambda,\lambda}\|_{X_{\lambda,\lambda}^{1-s,1-\theta}}^2\right)
\]
and due to \eqref{M} it remains to bound the right hand side by
\[
M + \| u\|_{L^2 H^{-s-\theta}}^2
\]
There is nothing to do for the $f_{\lambda,1}$ term. On the other hand
we can  bound
\[
\begin{split}
\|u_{\lambda,d}\|_{X^{1-s,1-\theta}_{\lambda,d}}^2 \lesssim &\
\lambda^{2(1-s)} d^{2(1-\theta)} \| u_{\lambda,d}\|_{L^2}^2
+ \lambda^{-2s} d^{-2\theta}  \| \tilde \Box_{g_{<
\sqrt\lambda}} u_{\lambda,d}\|_{L^2}^2 \\
= &\ \lambda^{2(1-s)} d^{2(1-\theta)} \| h_{\lambda,d} - h_{\lambda,2d}\|_{L^2}^2
+ \lambda^{-2s} d^{-2\theta}  \| f_{\lambda,d} - f_{\lambda,2d}\|_{L^2}^2
\\
\lesssim &\ \lambda^{-2s} d^{-2\theta}
 (\|f_{\lambda,d}\|_{L^2}^2+ \|f_{\lambda,2d}\|_{L^2}^2)
\\ &\  +    \lambda^{2(1-s)} d^{2(1-\theta)} (\|h_{\lambda,2d}\|_{L^2}^2
+\|h_{\lambda,d}\|_{L^2}^2)
\end{split}
\]
Finally, for the last term we have
\[
\begin{split}
 \| h_{\lambda,\lambda}\|_{X_{\lambda,\lambda}^{1-s,1-\theta}}^2
\lesssim  &\ \lambda^{2(1-s)} \lambda^{2(1-\theta)} \| h_{\lambda,\lambda}\|_{L^2}^2
+ \lambda^{-2s} \lambda^{-2\theta}  \| \tilde \Box_{g_{<
\sqrt\lambda}} h_{\lambda,\lambda}\|_{L^2}^2
\\ = &\ \lambda^{2(1-s)} \lambda^{2(1-\theta)} \| h_{\lambda,\lambda}\|_{L^2}^2
+ \lambda^{-2s} \lambda^{-2\theta}  \| S_\lambda u -  f_{\lambda,\lambda}\|_{L^2}^2
\\ \lesssim&\ \lambda^{2(1-s)} \lambda^{2(1-\theta)} \|
h_{\lambda,\lambda}\|_{L^2}^2 + \lambda^{-2s} \lambda^{-2\theta}  \| f_{\lambda,\lambda}\|_{L^2}^2
+ \| S_\lambda u\|_{L^2 H^{-s-\theta}}^2
\end{split}
\]
The proof is concluded.
\end{proof}

\section{Strichartz estimates and applications.}
\label{secse}

The Strichartz estimates for the variable coefficient wave
equation, as proved in \cite{MR1887639}, have the form:

\begin{theorem}(Tataru \cite{MR1887639}) Assume that the coefficients
  $g^{ij}$ of $\Box_g$ satisfy $\p^2 g^{ij} \in L^1L^\infty$. Then the
  solutions to the wave equation in $n+1$ dimensions satisfy the
  bounds
\begin{equation}
\| D^{\sigma} \nabla u\|_{L^p L^q} \lesssim \|u(0)\|_{H^1} + \|u_t(0)\|_{L^2} +
\|\Box_g u\|_{L^1 L^2}
\label{se}\end{equation}
where
\begin{equation}
\sigma = -\frac{n}2 + \frac{1}p +\frac{n}q, \qquad \frac{2}p +\frac{n-1}q
\leq \frac{n-1}2, \qquad 2 \leq p \leq \infty, \ \ 2 \leq q < \infty
\label{pq}\end{equation}
\end{theorem}

Applying this bound on an interval $I$ of size $\epsilon^2$ we
obtain by Cauchy-Schwartz
\[
\| D^{\sigma} \nabla u\|_{L^p(I; L^q)} \lesssim \frac{1}\epsilon \|u\|_{H^1(I
  \times \R^n)} + \epsilon \|\Box_g u\|_{L^2(I \times \R^n)}, \qquad
\epsilon \leq 1
\]
Summing up over small intervals this extends to intervals of arbitrary
lengths. Optimizing over $\epsilon$ yields
\begin{equation}
\| D^{\sigma} \nabla u\|_{L^p L^q}^2 \lesssim \| u\|_{H^1}^2 +
 \| u\|_{H^1} \|\Box_g u\|_{L^2}
\end{equation}
We want to apply this result to the functions $S_\l u_{\l,d}$ in
Definition~\ref{xst}.  By \eqref{strongst} we obtain

\begin{corollary}
a) Let $(\sigma,p,q)$ verifying
\[
\sigma = -\frac{n}2 + \frac{1}p +\frac{n}q, \qquad 2 \leq p \leq
\infty, \ \ 2 \leq q < \infty
\]
Then for $(\sigma,p,q)$ as in \eqref{pq} we have
\[
\| S_\l \nabla u\|_{L^p L^q}  \lesssim \lambda^{1-s-\sigma}
d^{\frac12-\theta} \|u\|_{X^{s,\theta}_{\l,d}}
\]
If additionally $\theta > \frac12$ then
\[
\| S_\l \nabla u\|_{L^p L^q}  \lesssim \lambda^{1-s-\sigma}
\|u\|_{X^{s,\theta}}
\]
b) If instead
\[
\frac{2}p +\frac{n-1}q \geq \frac{n-1}2
\]
then
\[
\| S_\l \nabla u\|_{L^p L^q}  \lesssim  \lambda^{1-
s-\sigma+\frac12(\frac{2}p +\frac{n-1}q - \frac{n-1}2)}
d^{\frac12-\theta-\frac12(\frac{2}p +\frac{n-1}q - \frac{n-1}2)}
\|u\|_{X^{s,\theta}_{\l,d}}
\]

\label{sob}
\end{corollary}

The interesting triplets of indices for $(\sigma,p,q)$ in $4+1$
dimensions are
\[
(0,\infty,2) \text{(energy)} \qquad (-\frac12,\frac{10}3,\frac{10}3)
\text{(Strichartz)}  \qquad (-\frac56,2,6) \text{(Pecher)}
\]
In addition,  we can also use the index $q=\infty$. Thus
we obtain the triplets
\[
(-2,\infty,\infty), \qquad (-\frac32,2,\infty)
\]
For the case when $\theta<\frac12$, we rely on the additional
triplets
\[
(-\frac16,2,3),\qquad (\frac14,4,2)
\]

For convenience we summarize the bounds we need for $\tilde
X^{s,\theta}_{\l,d}$:

\begin{corollary}
For $0 < \theta < 1$ we have
\[
\l^{s-1} \|S_\l \nabla u\|_{L^\infty L^2} + \l^{s-\frac52} \|S_\l \nabla u\|_{L^2
  L^\infty} + \l^{s-3} \| S_\l \nabla u\|_{L^\infty} \lesssim
d^{\frac12-\theta} \|u\|_{\tilde X^{s,\theta}_{\l,d}}
\]
\[
\l^{s-\frac{17}{12}} \|S_\l \nabla u\|_{L^2
  L^3} + \l^{s-1} \| S_\l \nabla u\|_{L^4 L^2} \lesssim
d^{\frac14-\theta} \|u\|_{\tilde X^{s,\theta}_{\l,d}}
\]
\label{d}\end{corollary}
The reason we include the gradient is to have also bounds for
$u_t$. Because of the frequency localization, if we drop the gradient
the same bounds hold with one less power of $\l$.

In our estimates later on we also need to work with $X^{s,b}$
functions which are concentrated into a smaller modulation range. For
this we introduce the additional norm
\[
\| u\|_{\tilde X^{s,\theta}_{\l,<d}}^2 = \inf \left\{
\sum_{h=1}^d \|u_h\|_{\tilde X^{s,\theta}_{\l,h}}^2;\ u = \sum_{h=1}^d u_h\right\}
\]
If $d = \l$ we simply write $\tilde X^{s,\theta}_{\l}$. A simple
argument leads to
\[
\|u\|^2_{X^{s,\theta}}\,=\,\inf \left \{ \sum_{\l=1}^\infty
 \|S_{\l} u_{\l}\|_{\tilde
X^{s,\theta}_{\l}}^2 ;\ u =\sum_{\l=1}^\infty S_{\l}
u_{\l}\right\}
\]
We also have

\begin{corollary}
a) Assume that $\theta > \frac12$. Then
\[
\l^{s-1} \|S_\l \nabla u\|_{L^\infty L^2} + \l^{s-\frac52} \|S_\l \nabla u\|_{L^2
  L^\infty} + \l^{s-3} \| S_\l \nabla u\|_{L^\infty} \lesssim  \|u\|_{\tX^{s,\theta}_{\l,<d}}
\]

b) Assume that $\theta < \frac12$. Then
\[
\l^{s-1} \|S_\l \nabla u\|_{L^\infty L^2} + \l^{s-\frac52} \|S_\l \nabla u\|_{L^2
  L^\infty} + \l^{s-3} \| S_\l \nabla u\|_{L^\infty} \lesssim
d^{\frac12-\theta} \|u\|_{\tilde X^{s,\theta}_{\l,<d}}
\]
\label{lessd}\end{corollary}

In preparation for proving bilinear estimates for the
$X^{s,\theta}$ spaces we first investigate which multiplications
leave the $\tX^{s,\theta}_{\l,d}$ space unchanged. For this we
define the algebras $M_d$, $M_{<d}$ with the norms
\[
\|f\|_{M_d} = \| f\|_{L^\infty} + d^{-1} \|f_t\|_{L^\infty}
 + d^{-\frac12} \|f_t\|_{L^2 L^\infty} + d^{-\frac32} \|f_{tt}\|_{L^2 L^\infty}
\]
\[
\|f\|_{M_{<d}} = \|f\|_{M_d} + d^{\frac12} \| f\|_{L^2  L^\infty}
\]
Then we have the multiplicative properties

\begin{lemma}
Assume that $f$ is localized at frequency $d \leq \l$.
Then we have
\begin{equation}
\|f S_\l u\|_{\tilde X^{s,\theta}_{\l,d}} \lesssim \|f\|_{M_d} \|u\|_{\tilde
  X^{s,\theta}_{\l,d}}
\end{equation}
respectively
\begin{equation}
\begin{split}
\|f S_\l u\|_{\tilde X^{s,\theta}_{\l,d}} & \lesssim \|f\|_{M_{<d}} \|u\|_{\tilde
  X^{s,\theta}_{\l,<d}}, \qquad \theta < \frac12 \\
\|f S_\l u\|_{\tilde X^{s,\theta}_{\l,d}} & \lesssim
d^{\theta-\frac12} \|f\|_{M_{<d}} \|u\|_{\tilde
  X^{s,\theta}_{\l,<d}}, \qquad \theta > \frac12
\end{split}
\end{equation}
\label{md}\end{lemma}

The proof is straightforward,  using Leibnitz's rule and the energy
estimate \eqref{litwo}. To bound functions in the $M_d$, respectively
$M_{<d}$ norms we use Corollary~\ref{lessd} with $d = \l$ to obtain:

\begin{lemma}
a) Assume that $\theta > \frac12$. Then
\[
\| S_\l u\|_{M_{<\l}} \leq \l^{2-s} \|u\|_{\tilde
X^{s,\theta}_{\l}}\qquad \| S_{<\l} u\|_{M_{\l}} \leq
\max\{1,\l^{2-s}\} \|u\|_{X^{s,\theta}}
\]
\[\| S_{<\l} u\|_{M_{<\l}} \leq
\max\{\l^ \frac12,\l^{2-s}\} \|u\|_{X^{s,\theta}}
\]

b) Assume that $\theta < \frac12$. Then
\[
\| S_\l u\|_{M_{<\l}} \leq \l^{\frac52-\theta-s} \|u\|_{\tilde
X^{s,\theta}_{\l}}\qquad \| S_{<\l} u\|_{M_{\l}} \leq
\max\{1,\l^{\frac52-\theta-s}\} \|u\|_{X^{s,\theta}}
\]
\[\| S_{<\l} u\|_{M_{<\l}} \leq
\max\{\l^ \frac12,\l^{\frac52-\theta-s}\} \|u\|_{X^{s,\theta}}
\]
\label{xmd}\end{lemma}

Using the above property we prove the algebra property \eqref{xx}
for the space $X$.

\begin{prop}
Assume that $s > 2$ and $\frac12 < \theta < s-\frac32$ .
Then $X^{s,\theta}$ is an algebra.
\label{xxx}\end{prop}

\begin{proof}
Let $u,v \in X^{s,\theta}$. For both we consider the decomposition
in Definition~\ref{xst},
\[
u  = \sum_{\l=1}^\infty \sum_{d=1}^{\l} S_\l u_{\l,d}, \qquad
v  = \sum_{\l=1}^\infty \sum_{d=1}^{\l} S_\l v_{\l,d},
\]
For the terms in the decomposition we  use the $\tilde{X}^{s,\theta}_{\l,d}$
norms,  as allowed by Corollary~\ref{tx}. We denote
\[
u_\l = \sum_{d=1}^{\l} u_{\l,d}, \qquad u_{\l,<d} =  \sum_{h=1}^{d} u_{\l,h}
\]
Then we write
\[
uv = \sum_{\mu =1}^\infty  S_{\mu} (uv) =  \sum_{\mu =1}^\infty
\sum_{\l_1 =1}^\infty  \sum_{\l_2 =1}^\infty S_{\mu} (S_{\l_1}u_{\l_1} S_{\l_2}v_{\l_2})
\]
There are two cases when the above summand
is nonzero, namely if $\l_1 \approx \l_2 \gtrsim \mu$ and if
$\max\{\l_1,\l_2\} \approx \mu$. We consider them separately.

{\bf Case 1}, $ \l_1,\l_2 \approx \l \gtrsim \mu$. In this case the
summability with respect to $\l$ is trivial, so it suffices to look at
the product $S_\l u_\l S_\l v_\l$ for fixed $\l$. This is localized at frequency
$\leq \l$.  Combining the $L^\infty L^2$ and the $L^2 L^\infty$ bounds
in Corollary~\ref{lessd} we obtain
\begin{equation}
\| S_\l u_\l S_\l v_\l\|_{L^2} +\l^{-1} \| \p_t (S_\l u_\l S_\l
v_\l)\|_{L^2} \lesssim   \l^{-2s+\frac32}
 \|u_\l\|_{\tilde X^{s,\theta}_{\l}} \|v_\l\|_{\tilde X^{s,\theta}_{\l}}
\label{lla}\end{equation}
Using the equation we can also bound the second time derivative,
\begin{equation}
\l^{-2} \| \p_t^2 (S_\l u_\l S_\l v_\l)\|_{L^2} \lesssim
\l^{-2s+\frac32}
 \|u_\l\|_{\tilde X^{s,\theta}_{\l}} \|v_\l\|_{\tilde X^{s,\theta}_{\l}}
\label{llb}\end{equation}
The three bounds above allow us to estimate for $\mu \leq \l$
\[
\|S_\l u_\l S_\l v_\l\|_{X^{s,\theta}_{\mu,\mu}} \lesssim
\mu^{s+\theta-2} \l^{-2s+\frac72}
 \|u_\l\|_{\tilde X^{s,\theta}_{\l}} \|v_\l\|_{\tilde X^{s,\theta}_{\l}}
\]
This suffices provided that $\theta < s-\frac32$, which is insured by
our hypothesis.

{\bf Case 2}. Here we consider products of the form $S_\mu v_\mu S_\l
u_\l $ where $\mu \ll \l$. Then the product is localized at frequency
$\l$. The summation with respect to $\l$ is trivial, but not the one
with respect to $\mu$.  We write
\[
S_\mu v_\mu S_\l u_\l  =  S_\mu v_{\mu} S_\l u_{\l,<\mu} +
\sum_{d = \mu}^\l  S_\mu v_{\mu} S_{\l} u_{\l,d}
\]
Using Lemma~\ref{md} and Lemma~\ref{xmd} we obtain
\[
\begin{split}
\|S_\l u_\l S_\mu v_\mu&\|_{\tilde X^{s,\theta}_{\l}}^2 \lesssim
\|S_\mu v_{\mu} S_\l u_{\l,<\mu}\|_{\tilde X^{s,\theta}_{\l,\mu}}^2
+  \sum_{d = \mu}^\l  \|S_\mu v_{\mu} S_{\l} u_{\l,d}\|_{\tilde X^{s,\theta}_{\l,d}}^2
\\ &\lesssim \mu^{2\theta-1}
\|S_\mu v_\mu\|^2_{M_{<\mu}}   \|u_{\l,<\mu}\|_{\tilde X^{s,\theta}_{\l,<\mu}}^2
+ \|S_\mu v_\mu\|^2_{M_{\mu}} \sum_{d=\mu}^\l  \|u_{\l,d}\|_{\tilde
  X^{s,\theta}_{\l,d}}^2
\\ &\lesssim  \mu^{2\theta-1}\|S_\mu
v_\mu\|^2_{M_{<\mu}}\sum_{d=1}^\l
 \|u_{\l,d}\|_{\tilde  X^{s,\theta}_{\l,d}}^2
\\ &\lesssim \mu^{3+2\theta-2s} \|v_\mu\|^2_{ \tX^{s,\theta}_{\mu}}\sum_{d=1}^\l
\|u_{\l,d}\|_{\tilde  X^{s,\theta}_{\l,d}}^2
\end{split}
\]
The summation with respect to $\mu$ is trivial since $\theta < s-\frac32$.

\end{proof}

We next prove \eqref{xy}.
\begin{prop}
Assume that $s > 2$ and $\frac12 < \theta < s-\frac32$ .
Then we have the multiplicative estimate
\[
X^{s,\theta} \cdot X^{s-1,\theta-1} \subset X^{s-1,\theta-1}
\]
\end{prop}

\begin{proof}
By duality this reduces to the multiplicative estimate
\[
X^{s,\theta} \cdot (X^{1-s,1-\theta}+L^2 H^{2-s-\theta})  \subset
X^{1-s,1-\theta} + L^2 H^{2-s-\theta}
\]
Since $s > 2$ we have the fixed time multiplication
\[
H^s  \cdot  H^{2-s-\theta} \subset H^{2-s-\theta}
\]
which implies the space-time bound
\[
L^\infty H^s \cdot L^2  H^{2-s-\theta} \subset L^2  H^{2-s-\theta}
\]
Due to the energy estimate for $X^{s,\theta}$ it remains to show that
\[
X^{s,\theta} \cdot X^{1-s,1-\theta} \subset
X^{1-s,1-\theta} + L^2 H^{2-s-\theta}
\]
We consider a product $S_\l u_\l S_\mu v_\mu$ which we decompose
as in the previous proof.  Because of the lack of symmetry we now
need to consider three cases.

{\bf Case 1}. Here we estimate $ S_\mu (S_\l u_\l S_\l v_\l)$ where $\mu
\lesssim \l$. By Corollary~\ref{lessd} we obtain
\[
\|S_\l u_\l S_\l v_\l\|_{L^2 L^\frac{3}{2}} \lesssim \|S_\l
u_\l\|_{L^4 L^3} \|S_\l v_\l\|_{L^4 L^3} \lesssim
\l^{\theta-\frac{2}{3}} \|u_\l\|_{\tilde X^{s,\theta}_{\l}}
\|v_\l\|_{\tilde
 X^{1-s,1-\theta}_{\l}}
\]

Using then Sobolev embeddings we obtain

\[
\|S_\mu(S_\l u_\l S_\l v_\l)\|_{L^2 H^{2-s-\theta}} \lesssim
\mu^{\frac 83 -s-\theta}\l^{\theta-\frac{2}{3}} \|u_\l\|_{\tilde
X^{s,\theta}_{\l}} \|v_\l\|_{\tilde
 X^{1-s,1-\theta}_{\l}}
\]

{\bf Case 2}. Here we bound $S_\mu u_\mu S_\l v_\l$, $\mu \ll \l$.
The product is localized at frequency $\l$, and the analysis is almost
identical to Case 2 in Proposition~\ref{xxx}.

{\bf Case 3}. Here we bound $S_\l u_\l S_\mu v_\mu$, $\mu \ll \l$.
The same argument applies, the only difference is that we gain some
extra $\mu/\l$ factors.

\end{proof}

We continue with the Moser estimates in \eqref{moser}, which follow from

\begin{prop}
Assume that $s > 2$ and $\frac12 < \theta < s-\frac32$ . Let $\Gamma$
be a smooth function.
Then
\[
\|\Gamma(u)\|_{X^{s,\theta}} \lesssim C(\|u\|_{L^\infty}) (1 +
\|u\|_{X^{s,\theta}}^5)
\]
\end{prop}
\begin{proof}
We write
\[
\Gamma(u)-\Gamma(v) = (u-v) f(u,v)
\]
and
\[
f(u,v) - f(x,y) = (u-x) g_1(u,v,x,y) + (v-y)g_2(u,v,x,y)
\]
where $f$, $g_1$ and $g_2$ are smooth functions.
Then we have
\[
\begin{split}
\Gamma(u) =& \Gamma(u_{1}) + \sum_{\l=1}^\infty \Gamma(u_{\leq 2\l}) -  \Gamma(u_{\leq \l})
\\ =& \Gamma(u_{1}) + \sum_{\l=1}^\infty u_{2\l} f(u_{\leq 2\l},u_{\leq \l})
\\ = &  \Gamma(u_{1}) + \sum_{\l=1}^\infty u_{2\l} [ f(u_{\leq 2},u_1) +
\sum_{\mu=2}^\l   ( f(u_{\leq 2\mu},u_{\leq \mu}) -f(u_{\leq
  \mu},u_{\leq \mu/2}))]
\\ =&  \Gamma(u_{1}) + \sum_{\l=1}^\infty u_{2\l} [ f(u_{\leq 2},u_1) +
\sum_{\mu=2}^\l  (u_{2\mu}\, g_1(u_{\leq 2\mu},u_{\leq
\mu},u_{\leq
  \mu/2}) \\ &+ u_\mu \, g_2(u_{\leq 2\mu},u_{\leq \mu},u_{\leq
  \mu/2}))]
\end{split}
\]
Hence we need  to bound
expressions of the form
\[
S_\l u_\l\, S_\mu v_\mu \, h(S_{<\mu} w), \qquad \mu \leq \l
\]
There are two different cases to consider:

{\bf Case 1}. $\mu \approx \l$. Then the product has the form
\[
S_\l u_\l\, S_\l v_\l \, h(S_{<\l} w)
\]
The first product is localized at frequency $\l$ and can be estimated
as in \eqref{lla}, \eqref{llb}.  For the nonlinear expression we
use Lemma~\ref{xmd} to obtain
\[
\| S_{<\l} w\|_{M_\l} \lesssim \|w\|_{X^{s,\theta}}
\]
On one hand by the chain rule we obtain
\begin{equation}
\| h(S_{<\l} w)\|_{M_\l} \lesssim C(\|w\|_{L^\infty}) (1 +
\|w\|_{X^{s,\theta}}^3) \label{mlx}\end{equation} On the other
hand because of the frequency localization we also have the
improved high frequency bound
\begin{equation}
\| \tS_\mu h(S_{<\l} w)\|_{M_\mu} \lesssim C(\|w\|_{L^\infty})
\left(\frac\l\mu\right)^N (1+ \|w\|_{X^{s,\theta}}^3), \qquad \mu
\gg \l \label{mlxmu}\end{equation} Taking this into account and
repeatedly using Leibnitz's rule we get
\[
\begin{split}
  & \| S_\l u_\l\, S_\l v_\l\, h(S_{<\l} w) \|_{X^{s,\theta}}^2 \\
  &\lesssim \sum_{\mu=1}^\infty \| S_\mu ( S_\l u_\l \,S_\l v_\l\,
  h(S_{<\l} w)) \|_{\tilde X^{s,\theta}_\mu}^2 \\ & \lesssim \sum_{\mu
    \lesssim \l} \| S_\l u_\l\, S_\l v_\l\, h(S_{<\l} w) \|_{\tilde
    X^{s,\theta}_\mu}^2 + \sum_{\mu \gg \l} \| S_\l u_\l\, S_\l
    v_\l\,
  \tS_\mu h(S_{<\l} w)) \|_{\tilde X^{s,\theta}_\mu}^2 \\ &
  \lesssim C(\|w\|_{L^\infty})
  \left (\sum_{\mu \lesssim \l} \mu^{2s+2\theta-4} \l^{-4s + 7} +
    \sum_{\mu \gg \l} \l^{2\theta + 3 - 2s} \left(\frac\l\mu\right)^N
  \right) \|u_\l\|^2_{\tilde X^{s,\theta}_\l} \|v_\l\|^2_{\tilde
    X^{s,\theta}_\l}  (1+\|w\|_{X^{s,\theta}}^6) \\ & \lesssim
  C(\|w\|_{L^\infty})\, \l^{2\theta + 3 - 2s} \, \|u_\l\|^2_{\tilde X^{s,\theta}_\l}
  \|v_\l\|^2_{\tilde X^{s,\theta}_\l} (1+\|w\|_{X^{s,\theta}}^6)
\end{split}
\]
This is trivially summable with respect to $\l$.

{\bf Case 2}. $\mu \ll \l$. Then the product has the form
\[
S_\l u_\l \,S_\mu v_\mu\, h(S_{<\mu} w) =
\]
\[
\begin{split}
&S_\l  u_{\l,<\mu}\, S_\mu v_\mu\, S_{<\mu} h(S_{<\mu} w) +
\sum_{\mu \leq
  d \ll \l} S_{\l}  u_{\l,<d}\, S_\mu v_\mu \, S_d h(S_{<\mu} w)
\\ & +  \sum_{\mu \leq d \ll \l} S_{\l}
u_{\l,d} \,S_\mu v_\mu \, S_{<d} h(S_{<\mu} w) + S_\l u_\l\, S_\mu
v_\mu\, S_{\l} h(S_{<\mu} w)
\\ & + \sum_{\nu \gg \l}  S_\l u_\l\, S_\mu v_\mu\, S_{\nu} h(S_{<\mu} w)
= f_1+f_2 +f_3 + f_4 + f_5
\end{split}
\]
For $f_1$ we use Lemma~\ref{md}, Lemma~\ref{xmd}  and \eqref{mlx}
to obtain
\[
\begin{split}
\|f_1\|_{\tilde X^{s,\theta}_{\l,\mu}} &\lesssim
 \|S_\l  u_{\l,<\mu}\, S_\mu v_\mu\|_{\tilde X^{s,\theta}_{\l,\mu}}
\| h(S_{<\mu} w)\|_{M_\mu} \\ &\lesssim \mu^{\theta-\frac12}
\|u_{\l,<\mu}\|_{\tilde X^{s,\theta}_{\l,<\mu}} \| S_\mu
v_\mu\|_{M_{<\mu}} \| h(S_{<\mu} w)\|_{M_\mu}
\\ &\lesssim C(\|w\|_{L^\infty}) \mu^{\theta+\frac32 -s}
\|u_{\l,<\mu}\|_{\tilde X^{s,\theta}_{\l,<\mu}} \|v_\mu\|_{\tilde
X^{s,\theta}_{\mu}} (1+\|w\|_{X^{s,\theta}}^3)
\end{split}
\]
The summation with respect to $\mu$ is trivial and the square
summability with respect to $\l$ is inherited from the first factor.

For $f_2$ we apply the same argument. There is a loss of a small power
of $(d/\mu)^\theta$ from the first product, but this is compensated by
the gain of arbitrary powers of $\mu/d$ due to \eqref{mlxmu}.
The same works for $f_3$ but there is no $(d/\mu)^\theta$ loss.
In the case of $f_4$ we need to worry about the $\l$ summability,
but the $(\mu/\l)^N$ gain in \eqref{mlxmu} settles this.
Finally, for $f_5$ there is a $(\mu/\nu)^N$ gain which cancels again
all the losses.

Summing up the pieces we obtain
\[
\begin{split}
&\|  S_\nu (S_\l u_\l \,S_\mu v_\mu\, h(S_{<\mu} w))
\|_{X^{s,\theta}}
\\ &\lesssim C(\|w\|_{L^\infty}) \nu^{s+\theta-2} \l^{-2s+\frac72} (\frac{\mu}{\l})^N \|u_\l\|_{\tilde X^{s,\theta}_\l }   \|v_\mu\|_{\tilde X^{s,\theta}_\mu }
(1+\|w\|_{X^{s,\theta}}^3)
\end{split}
\]
for $\nu \ll \l$,
\[
\|  S_\nu (S_\l u_\l\, S_\mu v_\mu \,h(S_{<\mu} w))
\|_{X^{s,\theta}} \lesssim C(\|w\|_{L^\infty}) \mu^{\theta+\frac32
-s} \|u_\l\|_{\tilde X^{s,\theta}_\l } \|v_\mu\|_{\tilde
X^{s,\theta}_\mu } (1+\|w\|_{X^{s,\theta}}^3)
\]
for $\nu \approx \l$, respectively
\[
\|  S_\nu (S_\l u_\l\, S_\mu v_\mu\, h(S_{<\mu} w))
\|_{X^{s,\theta}} \lesssim C(\|w\|_{L^\infty}) \mu^{\theta+\frac32
-s} (\frac{\mu}{\nu})^N \|u_\l\|_{\tilde X^{s,\theta}_\l }
\|v_\mu\|_{\tilde X^{s,\theta}_\mu } (1+\|w\|_{X^{s,\theta}}^3)
\]
for $\nu \gg \l$.

This concludes the proof of the proposition.

\end{proof}

Finally, we consider the bilinear estimate in \eqref{n}, which follows
from the next Proposition. Its proof cannot be completed using the
type of arguments we have employed so far. Instead, we contend
ourselves with reducing it to the trilinear estimate in
\eqref{finaltri}, to the proof of which we devote the rest of the
paper.

\begin{prop}
Assume that $s > \frac94$ and $\frac34 < \theta < s-\frac32$ .
Then we have the multiplicative estimate
\begin{equation}
\|\nabla u \nabla v\|_{X^{s-1,\theta-1}} \lesssim \|u\|_{X^{s,\theta}} \|v\|_{X^{s,\theta}}
\label{mainp} \end{equation}
\end{prop}

We begin our analysis with a simple observation, namely
that

\begin{lemma}
If $u \in X^{s,\theta}$ then $\nabla u \in \tilde X^{s-1,\theta}$ where
\[
\tilde X^{s-1,\theta} =  X^{s-1,\theta} + (L^2H^{s+\theta-1} \cap
H^1 H^{s+\theta-2}).
\]
\end{lemma}
\begin{proof}
We first consider spatial derivatives, for which we prove the better
bound
\[
\| \nabla_x u\|_{X^{s-1,\theta}} \lesssim \|u\|_{ X^{s,\theta}}
\]
By Definition~\ref{xst} and Corollary~\ref{tx} it suffices to show
that for functions $v$ localized at frequency $\lambda$ we have
\[
\| \nabla_x v\|_{X^{s-1,\theta}_{\lambda,d}} \lesssim
\|v\|_{\tilde X^{s,\theta}_{\lambda,d}}
\]
But this follows from the straightforward commutator bound
\begin{equation}
\| [  \Box_{g_{< \sqrt\lambda}},\nabla] v\|_{L^2} \lesssim
\lambda \|\nabla v\|_{L^2}
\label{ccv}\end{equation}
Here we recall that $g^{00}=1$, therefore every term in the commutator
has at least one spatial derivative.

Next we consider time derivatives, where it suffices to show that for
functions $v$ localized at frequency $\lambda$ we can write $v = v_1+
v_2$ where $v_1$, $v_2$ have the same frequency localization and
\begin{equation}
\| \partial_t v_1\|_{X^{s-1,\theta}_{\lambda,d}} +
\left(\frac{\lambda}{d}\right)^{1-\theta}
\|\partial_t v_2\|_{ (L^2H^{s+\theta-1} \cap
H^1 H^{s+\theta-2})} \lesssim
\|v\|_{\tilde X^{s,\theta}_{\lambda,d}}
\label{vud}\end{equation}
Roughly speaking $v_1$ accounts for the low modulation ($\lesssim
\lambda$) part of $v$ while $v_2$ accounts for the high modulation
part. We define $v_2$ as
\[
v_2 = (\Delta_{x,t})^{-1}  \Box_{g_{< \sqrt\lambda}} v
\]
This  satisfies the bound
\[
\|\nabla^2 v_2\|_{L^2} \lesssim
\|\Box_{g_{< \sqrt\lambda}} v\|_{L^2}
\]
which implies both the $v_2$ bound in \eqref{vud} and an $H^2$
bound for $ v_1$ which gives the correct $L^2$ bound for $\p_t
v_1$,
\[
\lambda^{s-1} d^\theta \| \partial_t v_1\|_{L^2} +
\left(\frac{\lambda}{d}\right)^{1-\theta} \|\partial_t v_2\|_{
  (L^2H^{s+\theta-1} \cap H^1 H^{s+\theta-2})} \lesssim \|v\|_{\tilde
  X^{s,\theta}_{\lambda,d}}
\]

 It remains to estimate
$ \Box_{g_{< \sqrt\lambda}} \partial_t v_1$. We have
\[
\|  \Box_{g_{< \sqrt\lambda}} \partial_t v_1 \|_{L^2} \leq  \| [
\Box_{g_{< \sqrt\lambda}},\partial_t] v_1\|_{L^2} + \|  \partial_t
\Box_{g_{< \sqrt\lambda}}  v_1 \|_{L^2}
\]
For the first term we use again \eqref{ccv}.
For the second  we compute
\[
\Box_{g_{< \sqrt\lambda}} v_1 = (-\Box_{g_{< \sqrt\lambda}}
+\Delta_{x,t}) (\Delta_{x,t})^{-1} \Box_{g_{< \sqrt\lambda}} v
\]
Since the difference $\Box_{g_{< \sqrt\lambda}} -\Delta_{x,t}$
contains no second order time derivatives this yields the bound
\[
\| \partial_t \Box_{g_{< \sqrt\lambda}} v_1\|_{L^2} \lesssim \lambda
\|\Box_{g_{< \sqrt\lambda}} v\|_{L^2}
\]
This allows us to conclude the proof of \eqref{vud} and therefore
the proof of the lemma.
\end{proof}

We now return to the estimate \eqref{mainp}.
Using the duality in \eqref{dual1}, \eqref{mainp} reduces to
\begin{equation}
\begin{split}
\left|\int uvw dx dt\right| \lesssim & \|u\|_{\tilde X^{s-1,\theta} }
\|v\|_{\tilde X^{s-1,\theta}}
\|w\|_{X^{1-s,1-\theta}+L^2 H^{2-s-\theta}}
\end{split}
\label{tri}\end{equation}
We do a trilinear Littlewood-Paley decomposition. Due to symmetry,
we need to consider two cases.

{\bf Case 1}. Here we consider high-high-low interactions and bound
\[
I = \int  S_\l u\, S_\l v\, S_\mu w \,dx dt, \qquad \mu \lesssim
\l
\]
We have
\[
|I| \lesssim \|S_\l u\|_{L^\infty L^2} \|S_\l v\|_{L^2 L^6}
\|S_\mu w\|_{L^2 L^3}
\]
which by the embeddings in Corollary~\ref{sob} give
\[
|I|  \lesssim \l^{\frac56-2s+2}  \mu^{s+\theta - \frac43}
\|u\|_{\tilde X^{s-1,\theta} }
\|v\|_{\tilde X^{s-1,\theta}}
\|w\|_{X^{1-s,1-\theta}+L^2 H^{2-s-\theta}}
\]
This suffices since both the exponent of $\l$ and the sum of the two
exponents are negative.

{\bf Case 2}. Here we consider high-low-high interactions and
seek to bound
\[
I = \int  S_\l u \,S_\mu v\, S_\l w \,dx dt, \qquad \mu \ll \l
\]
As a first simplification we dispense with the auxiliary $L^2$
norms. Begin with
\[
\begin{split}
|I| &\lesssim \|S_\l u\|_{L^2} \|S_\mu v\|_{L^\infty} \|S_\l w\|_{L^2}
\\ &\lesssim
\l^{s-1}\mu^{\frac34} \|S_\l u\|_{ L^2}\ \mu^{\frac94-s}
\|v\|_{\tilde X^{s-1,\theta}}\ \l^{1-s}\|S_\l w\|_{L^2}
\end{split}
\]
This allows us to dispense not only with the $L^2H^{s+\theta-1}$ part
of $u$, but also with its $X^{s-1,\theta}_{\l,>\mu}$ component.

If $v \in  L^2H^{s+\theta-1} \cap H^1 H^{s+\theta-2}$ then we
bound
\[
\begin{split}
|I| & \lesssim \|S_\l u\|_{L^\infty L^2} \|S_\mu v\|_{L^2
L^\infty} \|S_\l w\|_{L^2}
\\ & \lesssim \mu^{3-s-\theta} \|u\|_{\tilde X^{s-1,\theta} }
\ \|v\|_{L^2H^{s+\theta-1}}\ \l^{1-s}\|S_\l w\|_{L^2}
\end{split}
\]

Finally, if $w \in L^2 H^{2-s-\theta}$ then we can also estimate
\[
\begin{split}
|I| &\leq  \|S_\l u\|_{L^\infty L^2} \|S_\mu v\|_{L^2 L^\infty}
\|S_\l w\|_{L^2}
\\ & \lesssim \mu^{\frac32-s+\theta} \|u\|_{\tilde X^{s-1,\theta} }
\ \|v\|_{\tilde X^{s-1,\theta}}\ \l^{1-s} \mu^{1-\theta} \|S_\l
w\|_{L^2 }
\end{split}
\]
which suffices for both the $L^2 H^{2-s-\theta}$ and the
$X^{1-s,1-\theta}_{\l,> \mu}$ components of $w$.  Hence we have
reduced \eqref{tri} to the bound
\begin{equation}
|I| \lesssim \|S_\l u\|_{ X^{s-1,\theta}_{\l,<\mu }} \|S_\mu
v\|_{X^{s-1,\theta}} \|S_\l w\|_{X^{1-s,1-\theta}_{\l,<\mu}}
\qquad \mu \ll \l \label{tri1}
\end{equation}

Unfortunately we cannot fully prove this using Strichartz type
estimates. However, we can use scaling to simplify this further
and reduce it to
\begin{equation}
  \left|\int  S_\l u \,S_\mu v \,S_\l w dx dt\right| \lesssim \ln \mu
  \ \|u\|_{X^{0,1}_{\l,1}}
  \|v\|_{X^{\frac54,1}_{\mu,1}} \|w\|_{X^{0,\frac14}_{\l,d}} \qquad \mu
  \ll \l
\label{finaltri}
\end{equation}
For now we show that \eqref{finaltri} implies \eqref{tri1}. The
remaining sections of the paper are devoted to the proof of
\eqref{finaltri}.

After cancelling the powers of the high frequency the estimate
\eqref{tri1} follows after summation with respect to $1 \leq d_1, d_2,
d_3 \leq \mu$ from the bounds
\begin{equation}
\left|\int S_\l u \,S_\mu v \,S_\l w dx dt\right| \lesssim \ln \mu \
d_{min}^{\frac12} d_{mid}^\frac12 d_{max}^\frac14 \| u\|_{
X^{0,0}_{\l,d_1}} \|v\|_{ X^{\frac54,0}_{\mu,d_2}} \|
w\|_{X^{0,0}_{\l,d_3}}
\label{casea}\end{equation}
if $d_2 < d_{max}$, respectively
\begin{equation}
\left|\int S_\l u \,S_\mu v \,S_\l w dx dt\right| \lesssim \ln{\mu}\
d_{min}^{\frac12} d_{max}^\frac34 \|u\|_{  X^{0,0}_{\l,d_1}} \|
v\|_{ X^{\frac54,0}_{\mu,d_2}} \| w\|_{X^{0,0}_{\l,d_3}}
\label{caseb}\end{equation}
if $d_2 = d_{max}$.

To reduce all these cases to \eqref{finaltri} we use scaling combined
with a time decomposition argument. Precisely, for  $1 < d < \lambda$
we consider a smooth partition of unity in time with respect to
time intervals of length $d^{-1}$,
\[
1 = \sum \chi_d^j(t)
\]
Then a simple commutation argument shows that we can localize
the $\tilde X^{s,\theta}_{\lambda,d}$ norm to the $d^{-1}$ time
intervals while retaining square summability,
\begin{equation}
\| u\|_{ \tilde X^{s,\theta}_{\lambda,d}}^2 \approx
\sum_j \|  \chi_d^j u\|_{ \tilde X^{s,\theta}_{\lambda,d}}^2
\label{timedec}\end{equation}

We use such time decompositions in order to carry out
the following  three reduction steps:

(i) Reduction to $d_{min}=1$. By \eqref{timedec} all three norms are
square summable with respect to time intervals of length
$d_{min}^{-1}$. Hence it suffices to prove the bounds on
$d_{min}^{-1}$ time intervals. Rescaling such time intervals back to
time $1$ we arrive at the case $d_{min}=1$.  The regularity of the
coefficients improves after the rescaling, here and below.  Also we note that by
Duhamel's formula we can replace the factor corresponding to $d_{min}$
by a solution to the homogeneous equation.

(ii) Reduction to $d_{mid} =1$. By \eqref{timedec} the norms
corresponding to $d_{max}$ and $d_{mid}$ are square summable with
respect to time intervals of length $d_{mid}^{-1}$.  Hence it suffices
to prove the bounds on $d_{mid}^{-1}$ time intervals.  Rescaling such
time intervals back to time $1$ we arrive at the case $d_{mid}=1$.
Again by Duhamel's formula we also replace the factor corresponding to
$d_{mid}$ by a solution to the homogeneous equation.

(iii) Here we are in the case where two of the factors are solutions
for the homogeneous equation. In the case of \eqref{casea} the
remaining factor is at high frequency $\lambda$; then we use directly
\eqref{finaltri}.

In the case of \eqref{caseb} the remaining factor is at low frequency
$\mu$, so we need to prove that
\[
\left|\int S_\l u \,S_\mu v \,S_\l w dx dt\right| \lesssim \ln{\mu}\
 d^\frac34 \|u\|_{  X^{0,0}_{\l,1}} \|
v\|_{ X^{\frac54,0}_{\mu,d}} \| w\|_{X^{0,0}_{\l,1}}
\]
Partitioning the unit time into about $d$ time intervals of length
$d^{-1}$ this would follow from
\[
\left|\int \chi_d^i S_\l u \,S_\mu v \,S_\l w dx dt\right| \lesssim \ln{\mu}\
 d^\frac14 \|u\|_{  X^{0,0}_{\l,1}} \|
v\|_{ X^{\frac54,0}_{\mu,d}} \| w\|_{X^{0,0}_{\l,1}}
\]
Rescaling the small time intervals to unit size this becomes
exactly  \eqref{finaltri}.

\section{Half-waves and angular
localization operators}

We write the symbol for $\Box_g$,
\[
p(t,x,\tau,\xi) = \tau^2 - 2 g^{0j} \tau \xi_j - g^{ij} \xi_i
\xi_j
\]
in the form
\[
p(t,x,\tau,\xi) = (\tau + a^+(t,x,\xi)) (\tau + a^-(t,x,\xi))
\]
This leads to a decomposition of solutions to the wave equation
into two half-waves:

\begin{prop}(Geba-Tataru \cite{MR2153517})
  Let $u$ be a solution to the inhomogeneous equation \eqref{leq} for
  $\Box_g$. Then there is a representation
\[
\nabla u = u^+ + u^-
\]
where
\[
\begin{split}
\|u^+\|_{L^2} + \|(D_t+A^+(t,x,D)) u^+\|_{L^2} &+ \|u^-\|_{L^2} +
\|(D_t+A^-(t,x,D)) u^-\|_{L^2} \\ & \lesssim \|u\|_{H^1} +\|\Box_g
u\|_{L^2}
\end{split}
\]
\end{prop}

As a consequence, in \eqref{finaltri} we are allowed to replace
solutions to the $\Box_g$ equation by solutions to the $D_t+A^+$,
respectively $D_t+A^-$ equation. We also denote
\[
\aligned \|u\|_{X_\pm}\,=\,\|u\|_{L^2} &+ \|(D_t+A^\pm(t,x,D))
u\|_{L^2}\\
\|u\|_{X_{\pm,d}}\,=\,d^\frac14\|u\|_{L^2} &+
d^{-\frac34}\|(D_t+A^\pm(t,x,D)) u\|_{L^2}
\endaligned
\]

In order to facilitate the use of microlocal analysis tools it is
convenient to replace the symbols $a^\pm$ with mollified versions
$a^{\pm}_{<\mu}$ defined by
\[
a^{\pm}_{<\mu}(t,x,\xi) = S_{<\mu}(D_x) a(t,x,\xi)
\]

Given an angular scale $\alpha$ we consider the $\pm$ Hamilton flows
for $D_t+A^{\pm}_{<\alpha^{-1}}$.
\begin{equation}
\left\{ \begin{array}{c} \frac{d}{dt} x_t^\pm =  \p_\xi
a^\pm_{<\alpha^{-1}}(t,x_t^\pm,\xi_t^\pm) \cr \cr
 \frac{d}{dt} \xi_t^\pm =  -\p_x a^\pm_{<\alpha^{-1}}(t,x_t^\pm,\xi_t^\pm)
\end{array}\right . \qquad
 \left\{ \begin{array}{c} x_0^\pm=x \cr \cr \xi_0^\pm = \xi\end{array}
 \right .
\label{hf}
\end{equation}
These are bilipschitz flows, homogeneous with respect to the $\xi$
variable. The angular scale is relevant in that the Hamilton flow for
$D_t+A^{\pm}_{<\alpha^{-1}}$ serves as a good approximation to the
Hamilton flow for $D_t+A^{\pm}$ up to an $O(\alpha)$ angular
difference.

To characterize the higher regularity properties of these flows is
convenient to introduce (see \cite{TG}) a metric $g_\alpha$ in
the phase space, defined by
\[
ds^2 = |\xi|^{-4} (\xi d\xi)^2 + |\xi|^{-4} \alpha^{-2} (\xi \wedge
d\xi)^2 + \alpha^{-4} |\xi|^{-2} (\xi dx)^2 + |\xi|^{-2} \alpha^{-2}
(\xi \wedge dx)^2
\]
Then as in \cite{TG} we obtain
\begin{lemma}
The Hamilton flow maps $(x_t^\pm, \xi_t^\pm)$ are $g_\alpha$-smooth
canonical transformations.
\label{flowreg}\end{lemma}

Given a direction $\theta \in S^{n-1}$ at time $t=0$ we introduce the
size $\alpha$ sectors
\[
S_\alpha(\theta) = \{ \xi; \ \angle(\xi,\theta) < \alpha\}
\]
\[
\tilde S_\alpha(\theta) = \{ \xi; \ C \alpha < \angle(\xi,\theta) < 2
C\alpha\}
\]
where $C$ is a fixed large constant. The images
of $\R^n \times S_\alpha(\theta)$, respectively
 $\R^n \times \tilde S_\alpha(\theta)$  along the Hamilton
flow for $D_t+A^{\pm}_{<\alpha^{-1}}$ are denoted by $H_\alpha^\pm
S_\alpha(\theta)$, respectively $H_\alpha^\pm \tilde S_\alpha(\theta)$.

Let $\xi_\theta^\alpha =\xi_\theta^\alpha(x,t)$ be the Fourier
variable which is defined by the $D_t +A^{+}_{< \alpha^{-1}}$
Hamilton flow with initial data $\xi_\theta^\alpha(x,0)=\theta$
(i.e. $\xi_\theta^\alpha(x,t)=\xi^+_t(t)$ is the solution of the
flow \eqref{hf} with initial data $\xi^+_0=\xi$, for which
$x^+_t(t)=x$). This is well defined at least for a short time,
precisely for as long as caustics do not occur. From
Lemma~\ref{flowreg} one also sees that $\xi_\theta^\alpha$ is a
$g_\alpha$-smooth function of $x$.

We consider a maximal set $O_\alpha$ of $\alpha$-separated
directions and a partition of unity at time $0$
\[
1 =\sum_{\theta \in O_\alpha} \chi^{\pm,\alpha}_\theta(0,x,\xi)
\]
consisting of $0$-homogeneous symbols supported in
$S_\alpha(\theta)$ which are smooth on the corresponding scale.
Transporting these symbols along the $\pm$ Hamilton flows by
\[
 \chi^{\pm,\alpha}_\theta(0,x,\xi)\,=\,\chi^{\pm,\alpha}_\theta(t,x_t^\pm,\xi_t^\pm)
\]
produces a time dependent partition of unity
\begin{equation}
1 =\sum_{\theta \in O_\alpha} \chi^{\pm,\alpha}_\theta(t,x,\xi)
\label{linpart}
\end{equation}
so that the support of $\chi^{\pm,\alpha}_\theta(t,x,\xi)$ is
contained in $H_\alpha^\pm S_\alpha(\theta)$.

The regularity of these symbols is
easily obtained from the transport equations (see again \cite{TG}):

\begin{lemma}
The symbols $\chi^{\pm,\alpha}_\theta(t,x,\xi)$  belong to the
class $S(1,g_{\alpha})$\footnote{Throughout this paper we will use
the standard notation $S(m,g)$, while in \cite{TG} we used for
$S(1,g)$ the shorter one: $S(g)$.}. \label{wpsymbols}\end{lemma}

We use the above partition of unity in the phase space to
produce a corresponding pseudodifferential partition of unity.
Given a frequency $\lambda > \alpha^{-2}$ we define the symbols
\[
\chi^{\pm,\alpha}_{\theta,\lambda}(t,x,\xi) =  S_{<
  \lambda /8}(D_x) \chi^{\pm,\alpha}_{\theta}(t,x,\xi) \tilde s_\lambda(\xi)
\]
These are used in order to split general frequency localized waves
into square summable superpositions of directionally localized waves,
\[
S_\lambda u = \sum_{\theta \in
  O_\alpha}\chi^{\pm,\alpha}_{\theta,\lambda}(t,x,D) S_\lambda u
\]
This decomposition is closely related to a wave packet decomposition,
see \cite {MR1644105}, \cite{MR2178963}, \cite{MR2153517}, and
\cite{TG}. The difference is that here we skip the spatial
localization part since it brings no additional benefit. The above
localization at spatial frequencies less than $\lambda/8$ insures that
the output of the operators $\chi^{\pm,\alpha}_{\theta,\lambda}(t,x,D)
S_\lambda$ is still localized at frequency $\lambda$. This
localization is otherwise harmless:

\begin{lemma}
The symbols $\chi^{\pm,\alpha}_{\theta,\lambda} (t,x,\xi)$  belong to the class
$S(1,g_{\alpha})$. In addition, we have similar bounds for the Poisson
bracket
\begin{equation}
  \{ \tau+a^\pm_{<\alpha^{-1}} (t,x,\xi),
   \chi^{\pm,\alpha}_{\theta,\lambda} (t,x,\xi)\}
\in S(1,g_{\alpha})
\label{poib}\end{equation}
\label{wpsymbols}\end{lemma}

\begin{proof}
  The fact that $\chi^{\pm,\alpha}_{\theta,\lambda} (t,x,\xi) \in S(1,g_{\alpha})$ is straightforward since the
  multiplier $S_{<\lambda/8}$ is a mollifier on the $\lambda^{-1}$
  spatial scale, which is less that the spatial scale of the
  $g_\alpha$ balls.

Since $\chi^{\pm,\alpha}_{\theta}$ is transported along the
$a^\pm_{<\alpha^{-1}} (t,x,\xi)$ flow, the Poisson bracket is expressed in the form
\[
\{ a^\pm_{<\alpha^{-1}} (t,x,\xi), \tilde s_\lambda(\xi)\}
\chi^{\pm,\alpha}_{\theta,\lambda} (t,x,\xi) + \tilde s_\lambda(\xi) [
H_{a^\pm_{<\alpha^{-1}}},S_{<\lambda/8}(D_x)]
\chi^{\pm,\alpha}_{\theta}(t,x,\xi)
\]
Here $H_{a^\pm_{<\alpha^{-1}}}$ is the Hamiltonian operator
associated to the $a^\pm_{<\alpha^{-1}} (t,x,\xi)$ flow. It is
easy to see that the first term belongs to $S(1,g_\alpha)$,
therefore it remains to consider the commutator term. We have
\[
 [ H_{a^\pm_{<\alpha^{-1}}},S_{<\lambda/8}(D_x)] = [\partial_\xi
 a^\pm_{<\alpha^{-1}}, S_{<\lambda/8}(D_x)] \partial_x  - [\partial_x
 a^\pm_{<\alpha^{-1}}, S_{<\lambda/8}(D_x)] \partial_\xi
\]
The commutator of a scalar function $g$ with $S_{< \lambda/8}$
can be expressed as a rapidly convergent series of the form
\[
[g,S_{< \lambda/8}] = \lambda^{-1} \sum_j S_{< \lambda/8}^{1,j} \nabla g S_{<
  \lambda/8}^{2,j}
\]
where the multipliers $S_{< \lambda/8}^{1,j}$ and $S_{<
  \lambda/8}^{2,j}$ have the same properties as  $S_{<\lambda/8}$ and
decay rapidly with respect to $j$. Then the above commutator term
is expressed as
\[
 [ H_{a^\pm_{<\alpha^{-1}}},S_{<\lambda/8}(D_x)] =
\lambda^{-1} \sum_j S_{< \lambda/8}^{1,j} \left( \partial_x \partial_\xi
 a^\pm_{<\alpha^{-1}}  \partial_x  - \partial_x^2
 a^\pm_{<\alpha^{-1}}\partial_\xi \right)    S_{<\lambda/8}^{j,2}
\]
At this stage the effect of the mollifiers is negligible and we
can use the regularity properties of $a^\pm$ and
$\chi^{\pm,\alpha}_{\theta}$ to directly compute
\[
\tilde s_\lambda(\xi) [
H_{a^\pm_{<\alpha^{-1}}},S_{<\lambda/8}(D_x)]
\chi^{\pm,\alpha}_{\theta}(t,x,\xi)  \in S(\frac{1}{\alpha^2
\lambda},g_\alpha)
\]

\end{proof}

To better understand the phase space localization provided by
$\chi^{\pm,\alpha}_{\theta,\lambda} $ consider some point
$(x_0,t_0)$ and the corresponding center direction
$\xi_\theta^\alpha(x_0,t_0)$. A spatial unit $g_\alpha$ ball
$B_\theta^\alpha(x_0,t_0)$ centered at $(x_0,t_0)$ has
dimensions\footnote{Here $n$ stands for the space dimension}
$\alpha^2 \times \alpha^{n-1}$ with the long sides normal to
$\xi_\theta^\alpha(x_0,t_0)$. Within the ball
$B_\theta^\alpha(x_0,t_0)$, $\chi^{\pm,\alpha}_{\theta,\lambda}
\tilde S_\lambda$ localizes frequencies to a sector of angle
$\alpha$ centered at $\xi_\theta^\alpha(x_0,t_0)$. Thus the
frequencies are localized to a radial rectangle centered at
$\lambda \xi_\theta^\alpha(x_0,t_0)$ of size $\lambda \times
(\alpha \lambda)^{n-1}$. In this picture, angle $\alpha$ wave
packets correspond to a spatial localization on the scale of the
above ball $B_\theta^\alpha(x_0,t_0)$, constructed along a fixed
ray of the Hamilton flow.

The $g_\alpha$ metric restricted to frequency $\lambda$ is slowly
varying and temperate at frequencies \footnote{This corresponds to the
  classical wave packets which are localized on the scale of the
  uncertainty principle. Above this threshold we are dealing with
  generalized wave packets, which may have a more complex structure,
  see \cite{MR2153517} and \cite{TG}} $\lambda \geq\alpha^{-2}$, and
in our analysis we will always be above this threshold.  Hence there
is a good pseudodifferential calculus for operators with
$S(1,g_\alpha)$ symbols. The semiclassical parameter $h=h(\alpha,\lambda)$
in the $S(1,g_\alpha)$ calculus at frequency $\lambda$
is given by
\[
h(\alpha,\lambda) = (\alpha^2 \lambda)^{-1}
\]
The $S(1,g_\alpha)$ symbols at frequency
$\lambda$ satisfy the bounds
\begin{equation}
\left|(\xi_\theta^{\alpha} \partial_x)^\sigma
(\xi_\theta^{\alpha} \wedge \partial_x)^\beta
\partial_\xi^\nu  (\xi \partial_\xi)^\gamma q(t,x,\xi)\right|
\lesssim \alpha^{-2\sigma -|\beta|} (\alpha \lambda)^{-\nu}
\label{simbolga}\end{equation}
Due to the $L^2$ in time regularity of the second order derivatives of
the coefficients we also introduce the space of symbols $L^2
S(1,g_\alpha)$ which at frequency $\lambda$ satisfy
\begin{equation}
\left|(\xi_\theta^{\alpha} \partial_x)^\sigma
(\xi_\theta^{\alpha} \wedge \partial_x)^\beta
\partial_\xi^\nu  (\xi \partial_\xi)^\gamma q(t,x,\xi)\right|
\lesssim \alpha^{-2\sigma -|\beta|} (\alpha \lambda)^{-\nu} f(t)
\label{l2simbolga}\end{equation} for some $f \in L^2$. In all the
operators we consider here, the function $f$ is the same:
\begin{equation}
f(t) = M(\|\nabla^2 g(t)\|_{L^\infty})
\label{fM}\end{equation}
In some of our estimates we need to deal with two distinct scales at a
given frequency $\lambda$, namely the angular scale $\alpha$ and the
$\lambda^{\frac12}$ scale at which the coefficients are truncated.
Correspondingly we introduce additional symbol classes $C^k_\lambda
S(1,g_\alpha)$ of symbols $q$ localized at frequency $\lambda$ which
satisfy the $S(1,g_\alpha)$ bounds \eqref{simbolga} for
$\sigma+|\beta| \leq k$, respectively the weaker estimate
\begin{equation}
\left|(\xi_\theta^{\alpha} \partial_x)^\sigma
(\xi_\theta^{\alpha} \wedge \partial_x)^\beta
\partial_\xi^\nu  (\xi \partial_\xi)^\gamma q(t,x,\xi)\right|
\lesssim (\alpha^{-2\sigma -|\beta|}+\alpha^{-k}
\lambda^{\frac{\sigma+|\beta|}2}) (\alpha \lambda)^{-\nu}
\label{simbolgal}\end{equation} for $\sigma+|\beta| > k$.  There
is still a calculus for such symbols, since the above bounds are
stronger than the $S(1,g_{\lambda^{\frac12}})$ bounds.  The
related classes of symbols $L^2 C^k_\lambda S(1,g_\alpha)$ are
defined in a manner which is similar to \eqref{l2simbolga}.

Using the calculus for the above symbol classes one can prove that the
partition of unity in \eqref{linpart} yields an almost orthogonal
decomposition of functions, namely

\begin{prop}
Fix a frequency $\l$ and let $\alpha > \l^{-\frac12}$. Then for each
function $u$ which is localized at frequency $\lambda$ we have
\begin{equation}
\sum_{\theta \in O_\alpha} \| \chi^{\pm,\alpha}_{\theta,\lambda}(t,x,D)
u\|_{X_\pm}^2 \approx \| u\|_{X_\pm}^2
\label{linsumeq}\end{equation}
\label{linsum}\end{prop}

\begin{proof}
We only outline the proof, since this result is essentially contained
in \cite{MR2153517}.  There are two bounds to prove. The first
\begin{equation}
\sum_{\theta \in O_\alpha} \| \chi^{\pm,\alpha}_{\theta,\lambda}(t,x,D)
u\|_{L^2}^2 \approx \| u\|_{L^2} ^2
\label{aort}\end{equation}
follows from the almost orthogonality of the operators
$\chi^{\pm,\alpha}_{\theta,\lambda}(t,x,D)$. This in turn is due to
the almost disjoint supports\footnote{modulo tails which are rapidly
  decreasing on the $g_\alpha$ scale} of
$\chi^{\pm,\alpha}_{\theta,\lambda}$ and to the $S(1,g_{\alpha})$
calculus.

Consider now the second bound
\begin{equation}
\sum_{\theta \in O_\alpha} \| (D_t + A^\pm) \chi^{\pm,\alpha}_{\theta,\lambda}(t,x,D)
u\|_{L^2}^2 \approx \|(D_t + A^\pm)  u\|_{L^2}^2 +O(
\| u\|_{L^2}^2)
\label{aortc}\end{equation}
We first establish it with $A^\pm$ replaced by $A^\pm_{< \lambda^{\frac12}}$,
\begin{equation}
  \sum_{\theta \in O_\alpha} \| (D_t + A^\pm_{< \lambda^{\frac12}})
\chi^{\pm,\alpha}_{\theta,\lambda}(t,x,D)
  u\|_{L^2}^2 \approx \|(D_t + A^\pm_{< \lambda^{\frac12}})   u\|_{L^2}^2
+O( \| u\|_{L^2}^2)
\label{aortd}\end{equation}
Due to \eqref{aort} and the energy bound
\[
\| u\|_{L^\infty L^2}^2 \lesssim \|   u\|_{L^2}^2
+ \|   u\|_{L^2}\|(D_t + A^\pm_{< \lambda^{\frac12}})   u\|_{L^2}
\]
it suffices to prove the commutator estimate
\begin{equation}
\sum_{\theta \in O_\alpha}\| [ D_t+A^\pm_{<\lambda^{\frac12}},
\chi^{\pm,\alpha}_{\theta,\lambda}(t,x,D)] u\|_{L^2}^2 \lesssim \|
u\|_{L^\infty L^2}^2 \label{aorte}\end{equation} which we split
into two components.

For the low frequency part of the coefficients we use
a second order commutator
\begin{equation}
\sum_{\theta \in O_\alpha}\| [ D_t+A^\pm_{<\alpha^{-1}},
\chi^{\pm,\alpha}_{\theta,\lambda}(t,x,D)] u\|_{L^2}^2 \lesssim \|
u\|_{L^\infty L^2}^2
\label{lowcom}\end{equation}
For this it suffices to prove that
\begin{equation}
[ D_t+A^\pm_{<\alpha^{-1}}, \chi^{\pm,\alpha}_{\theta,\lambda}(t,x,D)]
\in OP L^2 S(1,g_\alpha)
\label{lowcoma}\end{equation}
The summation with respect to $\theta \in O_\alpha$ follows by
orthogonality since the symbols for the above commutators will retain
the rapid decay away from the support of
$\chi^{\pm,\alpha}_{\theta,\lambda}$. Here it is important that
\eqref{fM} applies uniformly.

Due to the Poisson bracket bound in \eqref{poib} it suffices to show
that
\[
[A^\pm_{<\alpha^{-1}},
\chi^{\pm,\alpha}_{\theta,\lambda}(t,x,D)] +
i \{a^\pm_{<  \alpha^{-1}}, \chi^{\pm,\alpha}_{\theta,\lambda}\} (t,x,D) \in OP L^2 S(1,g_\alpha)
\]
Due to the frequency localization of
$\chi^{\pm,\alpha}_{\theta,\lambda}$, only the values of
$a^\pm(t,x,\xi)$ in the region $|\xi| \approx \lambda$ can affect
the above operator. At this point it is no longer important that
$a^\pm_{<
  \alpha^{-1}}$ and $\chi^{\pm,\alpha}_{\theta,\lambda}$ are related.
We consider a rapidly convergent spherical harmonics expansion of
$a^{\pm}$,
\[
a^{\pm}(t,x,\xi) = \sum_j b_j (t,x) \phi_j(\xi)
\]
where $b_j$ have the same regularity as the coefficients $g^{ij}$
while $ \phi_j(\xi) $ are homogeneous of order $1$. It suffices to
consider a single term $b(t,x) \phi(\xi)$ in this expansion and show
that
\begin{equation}
[b_{<\alpha^{-1}} \phi(D),
\chi^{\pm,\alpha}_{\theta,\lambda}(t,x,D)] + i \{b_{<\alpha^{-1}} \phi
  , \chi^{\pm,\alpha}_{\theta,\lambda}\} (t,x,D) \in OPL^2S(1,g_\alpha)
\label{scom}\end{equation}
To see this we consider the commutators with $b$ and with $\phi$.
The commutator term with $b$ has the form
\[
C_b = ([b_{<\alpha^{-1}} ,
\chi^{\pm,\alpha}_{\theta,\lambda}(t,x,D)] +
i \{b_{<\alpha^{-1}},\chi^{\pm,\alpha}_{\theta,\lambda}\} (t,x,D))\phi(D)
\]
Since $\partial_x^2  b_{<\alpha^{-1}} \in L^2S(1,g_\alpha)$,   $\partial_\xi^2
\chi^{\pm,\alpha}_{\theta,\lambda} \in S(\alpha^{-2}
\lambda^{-2},g_\alpha)$ and $\phi \in S(\lambda,g_\alpha)$,
the $S(g_\alpha)$ calculus at frequency $\lambda$ yields
the better result $C_b \in OPL^2S(\alpha^{-2} \lambda^{-1},g_\alpha)$,
which is tight only when $\alpha= \lambda^{-\frac12}$.

The commutator term with $\phi$ has the form
\[
C_\phi = b_{<\alpha^{-1}} ([\phi(D),
\chi^{\pm,\alpha}_{\theta,\lambda}(t,x,D)] + i \{\phi
  , \chi^{\pm,\alpha}_{\theta,\lambda}\} (t,x,D))
\]
The $b_{<\alpha^{-1}}$ factor belongs to $S(1,g_\alpha)$ and can be
neglected. The argument for the remaining part is somewhat more
delicate since it hinges on the homogeneity of $\phi$. With $b = 1$
denote by $\xi$ the input frequency for $C_\phi$ and by $\eta$ the
output frequency. Due to the homogeneity of $\phi$ we have the representation
\begin{equation}
\phi(\eta) - \phi(\xi) = (\eta-\xi) \nabla \phi(\xi) +
\psi(\xi,\eta)(\xi \wedge (\xi-\eta))^2
\label{phixieta}\end{equation}
where $\psi$ is a smooth and homogeneous of order $-3$ matrix valued function. For
$|\xi|,|\eta| \approx \lambda$ we can separate variables in $\psi$ and
express it as a rapidly convergent series
\[
\psi(\xi,\eta) = \lambda^{-3}\sum_{j} \psi^1_j(\eta) \psi^2_j(\xi)
\]
This gives a representation for $C_\phi$ of the form
\[
C_\phi = \lambda^{-3}\sum_{j} \psi^1_j(D) ((\xi \wedge \partial_x)^2
\chi^{\pm,\alpha}_{\theta,\lambda})(t,x,D) \psi^2_j(D)
\]
Since $\chi^{\pm,\alpha}_{\theta,\lambda}(x,D) \in S(1,g_\alpha)$ we
obtain $(\xi \wedge \partial_x)^2 \chi^{\pm,\alpha}_{\theta,\lambda}
\in S(\lambda^2 \alpha^{-2},g_\alpha)$ which shows that $C_\phi \in
OPS(\alpha^{-2} \lambda^{-1},g_\alpha)$. This concludes the proof of
\eqref{scom} and thus the proof of \eqref{lowcoma}.

For the intermediate frequency
part of the coefficients we have a first order commutator estimate
\begin{equation}
\sum_{\theta \in O_\alpha}\| [ A^\pm_{\alpha^{-1} < \cdot <
  \lambda^{\frac12} }, \chi^{\pm,\alpha}_{\theta,\lambda}(t,x,D)]  u\|_{L^2}^2
 \lesssim \|   u\|_{L^\infty L^2}^2
\label{incom} \end{equation}
Together with \eqref{lowcom} this implies \eqref{aorte}.

This follows from first order commutator estimate
\begin{equation}
[ A^\pm_{\alpha^{-1} < \cdot < \lambda^{\frac12} },
\chi^{\pm,\alpha}_{\theta,\lambda}(t,x,D)] \in OPL^2 C^1_\lambda S(1,g_\alpha)
\label{incoma} \end{equation}
Indeed, for a scalar function $b$ we can estimate
\[
\alpha^{-2} \| b_{\alpha^{-1} < \cdot < \lambda^{\frac12} } \|_{L^2
  L^\infty}
+ \alpha^{-1} \| \partial_x b_{\alpha^{-1} < \cdot < \lambda^{\frac12} } \|_{L^2
  L^\infty}
\lesssim \|\partial^2 b\|_{L^2 L^\infty}
\]
Applied to the the symbol $a^\pm$ as a function of $x$ this shows that
\[
a^\pm_{\alpha^{-1} < \cdot <
  \lambda^{\frac12} }\in L^2C^2_\lambda S(\alpha^2 \lambda,g_\alpha)
\]
Since $\chi^{\pm,\alpha}_{\theta,\lambda} \in S(1,g_\alpha)$, the
estimate \eqref{incoma} follows by pdo calculus. The square summability
with respect to $\theta$ is again due to the almost disjoint supports
of the symbols $\chi^{\pm,\alpha}_\theta$.

It remains to pass from \eqref{aortd} to \eqref{aortc}. Due to
the energy bound
\[
\| u\|_{L^\infty L^2}^2 \lesssim \|   u\|_{L^2}^2
+ \|   u\|_{L^2}\|(D_t + A^\pm_{< \lambda^{\frac12}})   u\|_{L^2}
\]
this is a consequence of the estimate
\[
\|A^\pm_{> \lambda^{\frac12}} u\|_{L^2} \lesssim \| u\|_{L^\infty L^2}
\]
applied to both $u$ and $\chi^{\pm,\alpha}_{\theta,\lambda}(t,x,D)
u$. Using the spherical harmonics decomposition of the symbols
$a^{\pm}$ as above this reduces to the straightforward bound
\[
\| b_{> \lambda^{\frac12}} u\|_{L^2} \lesssim \lambda^{-1} \|\partial^2
b\|_{L^2 L^\infty}
\|u\|_{L^\infty L^2}
\]

\end{proof}

The frequency localization in $\chi^{\pm ,\alpha} _{\theta,\lambda}$
contributes to improved Strichartz type estimates above the critical
range of exponents.  Begin for instance with the endpoint $L^2 L^6$
Strichartz estimate
\begin{equation}
\|\chi^{\pm ,\alpha} _{\theta,\lambda}(t,x,D)  u\|_{L^2 L^6}
\lesssim \l^\frac56 \| u\|_{X_\pm} \label{specher}\end{equation}
Here the angular frequency localization plays no role.  However,
suppose we want to use Bernstein's inequality to replace this by
an $L^2 L^\infty$ estimate. Modulo rapidly decaying tails, within
each spatial $g_\alpha$ ball $B^\alpha_\theta(x_0,t_0)$ the
function $ \chi^{\pm ,\alpha} _{\theta,\lambda} (t,x,D) u$ is
frequency localized in a dyadic sector section of size $\lambda
\times (\alpha \lambda)^3$.  Then the constant in Bernstein's
inequality is
\[
[\lambda \times (\alpha \lambda)^3]^\frac16 = \lambda^\frac23 \alpha^\frac12
\]
Hence we obtain the better $L^2 L^\infty$ bound
\begin{equation}
  \|\chi^{\pm ,\alpha} _{\theta,\lambda}(t,x,D)  u\|_{L^2 L^\infty} \lesssim
  \alpha^\frac12 \lambda^\frac32 \| u\|_{X_\pm},
\qquad  \alpha > \lambda^{-\frac12}
\label{pecher}
\end{equation}
A simpler related uniform bound is derived directly from the energy
estimates,
\begin{equation}
  \|\chi^{\pm ,\alpha} _{\theta,\lambda}(t,x,D)  u\|_{L^\infty} \lesssim
  \alpha^\frac32 \lambda^2 \| u\|_{X_\pm},
\qquad  \alpha > \lambda^{-\frac12}
\label{pechera}
\end{equation}
A similar bound holds for the right hand side of the $\chi^{\pm
  ,\alpha} _{\theta,\lambda}(t,x,D)  u$ equation.
Indeed, for $u \in X_\pm$  we can write
\[
(D_t + A^\pm) \chi^{\pm ,\alpha}_{\theta,\lambda}(t,x,D) u =
(D_t + A^\pm_{<\l^\frac12}) \chi^{\pm ,\alpha}_{\theta,\lambda}(t,x,D)
u +  A^\pm_{>\l^\frac12} \chi^{\pm ,\alpha}_{\theta,\lambda}(t,x,D)
u
\]
The first term belongs to $L^2$ and has a similar
frequency localization as $\chi^{\pm ,\alpha}_{\theta,\lambda}(t,x,D)
u$. The second is estimated directly  using  \eqref{pecher}. This yields
\begin{equation}
  \|(D_t + A^\pm)  \chi^{\pm ,\alpha} _{\theta,\lambda}(t,x,D)  u\|_{L^2 L^\infty} \lesssim
  \alpha^\frac32 \lambda^2 \| u\|_{X_\pm},
\qquad  \alpha > \lambda^{-\frac12}
\label{pecherb}
\end{equation}

Another way of taking advantage of the angular localization is in
corresponding bounds for derivatives. Consider the differentiation
operators $\xi_\theta^\alpha \wedge D$ whose symbol vanishes in the
$\xi_\theta^\alpha$ direction. Then in the support of $\chi^{\pm ,\alpha} _{\theta,\lambda}$
these symbols have size $\alpha \lambda$. Hence from \eqref{specher}
we also obtain
\begin{equation}
\|(\xi_\theta^\alpha \wedge D)\chi^{\pm ,\alpha} _{\theta,\lambda}(t,x,D) u\|_{L^2
  L^6} \lesssim (\alpha \l) \l^\frac56 \| u\|_{X_\pm}
\label{spechera}\end{equation}
We can argue in the same way for the energy estimates or for the $L^2
L^\infty$ bound in \eqref{pecher}. For convenience we collect
several such bounds in a single norm,
\[
\begin{split}
  \| v\|_{X_{\pm}^{\lambda, \alpha, \theta}} =&\ \|v\|_{X_{\pm}} +\|
  v\|_{L^\infty L^2} + \lambda^{-\frac56} \| v\|_{L^2 L^6}+
  \alpha^{-\frac12} \lambda^{-\frac32} \| v\|_{L^2 L^\infty}
+ \alpha^{-\frac32} \lambda^{-2} \|v\|_{L^\infty}
\\ &\ + \alpha^{-\frac32} \lambda^{-2} \| (D_t + A^\pm) v\|_{L^2 L^\infty}
 +   (\alpha
  \lambda)^{-1}\| (\xi_\theta^\alpha \wedge D) v\|_{L^\infty
    L^2}  \\ &\ + (\alpha
  \lambda)^{-1} (\lambda^{-\frac56} \| (\xi_\theta^\alpha \wedge D) v\|_{L^2 L^6} +
  \alpha^{-\frac12} \lambda^{-\frac32} \| (\xi_\theta^\alpha \wedge D)
  v\|_{L^2 L^\infty})
\end{split}
\]
and use it to state a corresponding version of \eqref{linsumeq},
\begin{equation}
\sum_{\theta \in O_\alpha} \| \chi^{\pm,\alpha}_{\theta,\lambda}(t,x,D)
u\|_{X_\pm^{\lambda,\alpha,\theta}}^2 \approx \|\tilde S_\l u\|_{X_\pm}^2
\label{limsumeqa}\end{equation}

We want to replace the partition of unity in \eqref{linpart} first
with a bilinear one and next with a trilinear one. Given two
frequencies $\mu < \lambda$, we denote $\alpha_\mu= \mu^{-\frac12}$
and introduce a corresponding bilinear partition of unity which is
useful when estimating the frequency $\mu$ output of the product
of two frequency $\lambda$ waves.  The main contribution
corresponds to opposite frequencies $\xi$ and $\eta$, therefore
we organize the following decomposition based on the dyadic
angle $\alpha_\mu \leq \alpha \leq 1$ between $\xi$ and $-\eta$.
Precisely, by superimposing the $\alpha$ angular decompositions
for $\alpha$ in the above range we obtain
\[
\begin{split}
& \tilde s_\lambda(\xi) \tilde s_\lambda(\eta)=
 \\ & \sum_{\theta_1,\theta_2\in
O_{\alpha_\mu}}^{|\theta_1 + \theta_2| \leq 2C
  \alpha_\mu} \ \sum_{\theta_3,\theta_4 \in
  O_{2\alpha_\mu}}^{|\theta_3 +\theta_4| \leq 4C
   \alpha_\mu}
\chi^{\pm, \alpha_\mu}_{\theta_1,\lambda}(t,x,\xi)
\chi^{\mp, \alpha_\mu}_{\theta_2,\lambda}(t,x,\eta)
\chi^{\pm,2 \alpha_\mu}_{\theta_3,\lambda}(t,x,\xi)
\chi^{\mp,2 \alpha_\mu}_{\theta_4,\lambda}(t,x,\eta)
\\ &+ \!\!
\sum_{\alpha=\alpha_\mu}^1 \!\!
\sum_{\theta_1,\theta_2 \in O_\alpha}^{C \alpha
  \leq |\theta_1+\theta_2| \leq 2 C
  \alpha}
\ \sum_{\theta_3,\theta_4 \in
O_{2\alpha}}^{|\theta_3+\theta_4| \leq 4C
  \alpha} \!\!
\chi^{\pm,\alpha}_{\theta_1,\lambda}(t,x,\xi)
\chi^{\mp,\alpha}_{\theta_2,\lambda}(t,x,\eta)
\chi^{\pm,2\alpha}_{\theta_3,\lambda}(t,x,\xi)
\chi^{\mp,2\alpha}_{\theta_4,\lambda}(t,x,\eta)
\end{split}
\]
To shorten this expression we redenote factors and
harmlessly simplify the summation notations to
\begin{equation}
1 =  \sum_{\theta \in O_{\alpha_\mu}}
  \phi^{\pm,\alpha_\mu}_{\theta,\lambda}(t,x,\xi) \phi^{\mp,\alpha_\mu}_{-\theta,\lambda}(t,x,\eta)
+ \sum_{\alpha=\alpha_\mu}^1 \  \sum_{\theta \in O_\alpha}
 \phi^{\pm,\alpha}_{\theta,\lambda}(t,x,\xi)
\tphi^{\mp,\alpha}_{-\theta,\lambda}(t,x,\eta)
\label{bilinsum}\end{equation}
where the tilde in $ \tphi^{\pm,\alpha}_{\theta,\lambda}$ indicates an
$O(C\alpha)$ angular separation from $\theta$.  The symbols
$\phi^{\pm,\alpha}_{\theta,\lambda}$, respectively
$\tphi^{\pm,\alpha}_{\theta,\lambda}$ retain the same properties as
$\chi^{\pm,\alpha}_{\theta,\lambda}$, namely
\begin{equation}
\phi^{\pm,\alpha}_{\theta,\lambda} \in
S(1,g_\alpha),
\qquad   \{ \tau+a^\pm_{<\alpha^{-1}} (t,x,\xi),
   \phi^{\pm,\alpha}_{\theta,\lambda} (t,x,\xi)\}
\in S(1,g_{\alpha})
\end{equation}
and  the same for $\tphi^{\pm,\alpha}_{\theta,\lambda}$.
In particular the counterpart of \eqref{limsumeqa} is still valid,
\begin{equation}
\sum_{\theta \in O_\alpha} \| \phi^{\pm,\alpha}_{\theta,\lambda}(t,x,D)
u\|_{X_\pm^{\lambda,\alpha,\theta}}^2 +
 \| \tphi^{\pm,\alpha}_{\theta,\lambda}(t,x,D)
u\|_{X_\pm^{\lambda,\alpha,\theta}}^2\approx \|\tilde S_\l u\|_{X_\pm}^2
\label{limsumeqb}\end{equation}

Finally, we arrive at the main trilinear symbol decomposition.
Its aim is to achieve a simultaneous angular decomposition
in trilinear expressions of the form
\[
\int u v w dx dt
\]
We denote the three corresponding frequencies by $\xi, \eta$ and
$\zeta$. We assume that each of the factors has a dyadic frequency
localization,
\[
|\xi| \approx |\eta| \approx \lambda, \qquad |\zeta| \approx \mu,
\qquad 1 \ll \mu \leq \lambda
\]
If the trilinear decomposition were translation invariant then only
its structure on the diagonal $\xi+\eta+\zeta=0$ is relevant. However,
in our case we are working with variable coefficient operators
therefore a neighborhood of the diagonal is relevant. The size of
this neighborhood is determined by the spatial regularity of the
symbols via the uncertainty principle.

Corresponding to  the first term in \eqref{bilinsum}
we consider a decomposition in $\zeta$ with respect to the dyadic
angle between $\zeta$ and $\theta$,
\[
\tilde s_\mu(\zeta) = \phi^{\pm,\alpha_\mu}_{\theta,\mu} (t,x,\zeta)
+ \sum_{\alpha > \alpha_\mu}   \tphi^{\pm,\alpha}_{\theta,\mu}
(t,x,\zeta)
\]
To understand the $\zeta$ decomposition corresponding to the
second term in \eqref{bilinsum} we first identify the location of the
diagonal $\xi+\eta+\zeta = 0$. Given the above dyadic localization of
$\xi,\eta$ and $\zeta$, if the angle between $\xi$ and $-\eta$
is of order $\alpha$, then the angle between $\xi$ and $\pm \zeta$
must be of order $\alpha \lambda \mu^{-1}$ which is larger than
$\alpha$. Thus the interesting angular separation threshold for $\zeta$
is $\alpha \lambda \mu^{-1}$. It would appear that there are
two cases to consider, namely when the angle between $\xi$ and $\zeta$
is small, and when the angle between $-\xi$ and $\zeta$
is small. However, due to our choice of the $\pm$ signs corresponding
to $\xi$, $\eta$ and $\zeta$, the latter case leads to nonresonant
wave interactions and loses its relevance. Hence, the significant
dyadic parameter here is the angle between $\xi$ and $\zeta$,
and the $\zeta$ decomposition has the form
\[
\tilde s_\mu(\zeta) = \phi^{\pm,\alpha \mu^{-1}
\lambda}_{\theta,\mu} (t,x,\zeta) +  \tphi^{\pm,\alpha \mu^{-1}
\lambda}_{\theta,\mu} (t,x,\zeta) + \sum_{\beta
>\alpha \mu^{-1}
\lambda }   \tphi^{\pm,\beta}_{\theta,\mu} (t,x,\zeta)
\]
Then the full trilinear decomposition has the form
\begin{equation}
\begin{split}
\tilde s_\l(\xi) \tilde s_\lambda(\eta) \tilde s_\mu(\zeta)
  = & \sum_{\theta \in O_{\alpha_\mu}}
\phi^{\pm,\alpha_\mu}_{\theta,\lambda}(t,x,\xi)
 \phi^{\mp,\alpha_\mu}_{-\theta,\lambda}(t,x,\eta)
 \phi^{\pm,\alpha_\mu}_{\theta,\mu}(t,x,\zeta)
\\ +& \sum_{\theta \in O_{\alpha_\mu}}  \phi^{\pm,\alpha_\mu}_{\theta,\lambda}(t,x,\xi)
 \phi^{\mp,\alpha_\mu}_{-\theta,\lambda}(t,x,\eta)
\sum_{\alpha > \alpha_\mu}
\tilde \phi^{\pm,\alpha}_{\theta,\mu}(t,x,\zeta)
\\ +& \sum_{ \alpha > \alpha_\mu }\sum_{\theta \in
  O_{\alpha}}
\phi^{\pm,\alpha}_{\theta,\lambda}(t,x,\xi) \tilde
\phi^{\mp,\alpha}_{-\theta,\lambda}(t,x,\eta) \tilde \phi^{\pm,\alpha
\mu^{-1} \lambda}_{\theta,\mu}(t,x,\zeta)
\\ +& \sum_{ \alpha > \alpha_\mu }\sum_{\theta \in
  O_{\alpha}}
\phi^{\pm,\alpha}_{\theta,\lambda}(t,x,\xi) \tilde
\phi^{\mp,\alpha}_{-\theta,\lambda}(t,x,\eta) \phi^{\pm,\alpha \mu^{-1}
\lambda}_{\theta,\mu}(t,x,\zeta)
\\ +& \sum_{ \alpha > \alpha_\mu }\sum_{\theta \in
  O_{\alpha}}
\phi^{\pm,\alpha}_{\theta,\lambda}(t,x,\xi) \tilde
\phi^{\mp,\alpha}_{-\theta,\lambda}(t,x,\eta) \sum_{\beta > \alpha \mu^{-1}
\lambda} \tilde\phi^{\pm,\beta}_{\theta,\mu}(t,x,\zeta)
\end{split}
\label{trilindec}\end{equation}
In the above sum the first three terms are the main ones, as they
account for the behavior near the diagonal. The remaining terms
have off diagonal support, and their contribution to trilinear forms
as above is negligible.

\section{Proof of the trilinear
estimate~\eqref{finaltri}}

As noted in the previous section, we can replace the spaces
$X^{s,\theta}_{\l,d}$ in \eqref{finaltri} with the $X_{\pm}$ spaces.
Hence we restate  \eqref{finaltri} in the form

\begin{prop}
For any choice of the $\pm$ signs and $1 < d < \mu \ll \lambda$ we have
\begin{equation}
\left|\int S_\l u\, S_\l v\, S_\mu w  dx dt \right | \lesssim \ln \mu \cdot
\mu^\frac54\|S_\lambda u\|_{X_{\pm}} \|S_\lambda v\|_{X_{\pm,d}}
\|S_\mu w\|_{X_\pm}
\label{lastone}\end{equation}
\end{prop}

\begin{proof}
We begin with several simple observations. First, by localizing to
a fixed smaller space-time scale and rescaling back to unit scale we can
insure that the coefficients $g^{ij}$ vary slowly inside a unit
cube,
\[
|\nabla_{x,t} g^{ij}| \ll 1
\]
This in turn insures that the Fourier variable does not vary much
along the Hamilton flow,
\[
|\xi_\theta^\alpha - \theta| \ll 1
\]

We can also localize all factors in frequency to angular regions of
small size, say $< \frac{1}{20}$. The corresponding localization
multipliers are easily seen to be bounded in $X_{\pm}$ and
$X_{\pm,d}$.

If the first two $\pm$ signs are identical then
the product $S_\l u\, S_\l v$ is concentrated at a time frequency of
the order of $\l$ which makes it almost orthogonal to $S_\mu w$,
hence the estimate above is much easier. Therefore without any
restriction in generality we fix the first sign to $+$ and the
second one to $-$. Even though the problem is not symmetric with
respect to the first two factors, the sign in the third factor
plays no role whatsoever, so we fix it to $+$. We denote
\[
a(t,x,\xi)\,=\,a^+(t,x,\xi)
\]
Then
\[
a^-(t,x,\xi)\,=\,-a(t,x,-\xi)
\]
We note that, for the purpose of the above estimates, in the
definition of $X_{\pm}$ at frequency $\lambda$ we can replace the
symbols $a(x,\xi)$ with their regularized versions,
namely $a_{<\l^\frac12}(x,\xi)$.

To keep the number of parameters small we first present the argument
in the case when $d=1$. Once this is done, we show what changes
are necessary for $d > 1$.

{\bf Case 1}: $d = 1$. Corresponding to the trilinear symbol
decomposition \eqref{trilindec} of the identity we consider the
corresponding pseudodifferential decomposition of the trilinear
expression in \eqref{lastone}. The we estimate each of
the five terms. We remark that, since $S_\l u$, $S_\l v$ and $S_\mu w$
are frequency localized in a small  angle, so are all the factors
in \eqref{lastone}.

{\bf Case 1, term I:}
\[
I = \sum_{\theta \in O_{\alpha_\mu}} \int
\phi^{+,\alpha_\mu}_{\theta,\lambda}(t,x,D) S_\l u\
 \phi^{-,\alpha_\mu}_{-\theta,\lambda}(t,x,D) S_\l v\
 \phi^{+,\alpha_\mu}_{\theta,\mu}(t,x,D) S_\mu w \,dx dt
\]
We use the energy estimate for the first two factors and the
$L^2 L^\infty$ bound for the third to obtain
\[
|I| \lesssim \mu^{\frac54} \|
\phi^{+,\alpha_\mu}_{\theta,\lambda}(x,D) S_\l
u\|_{X^{\lambda,\alpha_\mu,\theta}_{+}} \|
\phi^{-,\alpha_\mu}_{-\theta,\lambda}(x,D) S_\l
v\|_{X^{\lambda,\alpha_\mu,\theta}_{-}} \|
\phi^{+,\alpha_\mu}_{\theta,\mu}(t,x,D) S_\mu w
\|_{X^{\mu,\alpha_\mu,\theta}_{+}}
\]
 The summation with
respect to $\theta$ is straightforward due to \eqref{limsumeqb}.

{\bf Case 1, term II:}
This is the most difficult term,
\[
II = \sum_{\theta \in O_{\alpha_\mu}} \int
\phi^{+,\alpha_\mu}_{\theta,\lambda}(t,x,D) S_\l u\
 \phi^{-,\alpha_\mu}_{-\theta,\lambda}(t,x,D)S_\l  v
\sum_{\alpha > \alpha_\mu} \tilde
\phi^{+,\alpha}_{\theta,\mu}(t,x,D) S_\mu w \,dx dt
\]

The summation with respect to $\theta$ is easily done using
\eqref{limsumeqb}.  Hence, in what follows, we fix $\theta$ and
redenote
\[
u_\theta = \phi^{+,\alpha_\mu}_{\theta,\lambda}(t,x,D) S_\l u, \qquad
v_\theta = \phi^{-,\alpha_\mu}_{-\theta,\lambda}(t,x,D) S_\l v, \qquad
w_{\theta}^\alpha = \tphi^{+,\alpha}_{\theta,\mu}(t,x,D)S_\mu  w
\]
The factors $u_\theta$ and $v_\theta$ are frequency localized
in small angles around $\theta$, respectively $-\theta$;  $w_{\theta}^\alpha$
has a similar localization around $\pm \theta$ provided that $\alpha
\ll 1$.

We denote by $\ta_{<\mu^\frac12} (t,x,\xi)$ the linearization of
$a_{<\mu^\frac12}(t,x,\xi)$ with respect to $\xi$ around
$\xi = \xi_\theta^{\alpha_\mu}(t,x)$. Since $a_{<\mu^\frac12}(t,x,\xi)$ is a
homogeneous symbol of order $1$, we have
\[
\ta_{<\mu^\frac12} (t,x,\xi) = \xi \partial_\xi
a_{<\mu^\frac12}(t,x,\xi_\theta^{\alpha_\mu})
\]
Consider now the difference
\[
e =  a_{<\mu^\frac12}- \ta_{<\mu^\frac12}
\]
It vanishes of second order on the half line $\R^+ \xi_\theta$.  Due
to the uniform (nonradial) convexity of the characteristic cone $\{
\tau + a_{<\mu^\frac12}(t,x,\xi)=0\}$, it follows that $e$ is nonzero
when $\xi$ is not collinear with $\xi_\theta^{\alpha_\mu}$.
Precisely, we can estimate it in terms of the angle
$\angle(\xi,\xi_\theta^{\alpha_\mu})$ as
\[
e(t,x,\xi) \approx  |\xi|
|\angle(\xi,\xi_\theta^{\alpha_\mu})|^2
\]
In particular in the support of the symbol
$\tphi^{+,\alpha}_{\theta,\mu}$ the above angle has size $\alpha$
and the frequency has size $\mu$.
Hence\footnote{here we switch to the letter $\zeta$ for the frequency,
as the following analysis refers to the  region at low frequency $\mu$
corresponding to the last factor $w$ in the trilinear form.}
\[
 e(t,x,\zeta)
\approx \alpha^2 \mu, \qquad (t,x,\zeta) \in \text{supp }
 \tphi^{+,\alpha}_{\theta,\mu}
\]
Here it may help to think of the constant coefficient case where
$\xi_\theta^{\alpha_\mu}=\theta$, while $a-\ta = |\xi| - \xi \theta$.
We introduce a local inverse for $ e(t,x,\zeta)$ in the support of
$\tphi^{+,\alpha}_{\theta,\mu}$, namely
\[
l(t,x,\zeta) = \tilde\tphi^{+,\alpha}_{\theta,\mu}(t,x,\zeta)
e^{-1}(t,x,\zeta)
\]
The cutoff symbol $ \tilde\tphi^{+,\alpha}_{\theta,\mu}$ is similar to
$\tphi^{+,\alpha}_{\theta,\mu}$ but has a slightly larger support and
equals $1$ in a neighbourhood of the support of
$\tphi^{+,\alpha}_{\theta,\mu}$.

As defined, the operator $L(t,x,D)$ is not localized at frequency
$\mu$. To remedy this we truncate its output in frequency and set
\[
\tilde L = \tilde S_\mu(D) L(t,x,D)
\]
 The properties of the operator
$\tilde L$ are summarized in the following

\begin{lemma}
The operator $\tilde L$ satisfies the following estimates:

a) fixed time $L^p$ mapping properties:
\[
\| \tilde L \|_{L^p \to L^p} \lesssim \alpha^{-2} \mu^{-1}, \qquad 1 \leq p \leq \infty
\]
b) fixed time approximate inverse of $A(t,x,D) - \tA(t,x,D)$:
\[
\|(A(t,x,D)- \tA(t,x,D)) \tilde L - \tilde\tphi(t,x,D)\|_{L^p
  \to L^p} \lesssim \mu^{-\frac12} +\alpha^{-2} \mu^{-1}, \quad 1 \leq
p \leq \infty
\]
c) space-time $X_+$ mapping properties:
\[
\| \tilde L  \|_{X_{+} \to X_+} \lesssim \alpha^{-2} \mu^{-1}
\]

\label{L}\end{lemma}
\begin{proof}
  We first compute the regularity of the symbol $e(t,x,\zeta)$ within
  the support of $l$. With respect to $\xi$ this is smooth and
  homogeneous, therefore we only have to keep track of the order of
  vanishing when $\xi$ is in the $\xi_\theta^{\alpha_\mu}$ direction.
  With respect to $x$ there is the dependence coming from the symbol
  $a$, as well as the dependence due to the $\xi_\theta^{\alpha_\mu}$
  direction occuring in the linearization. Since $a$ is Lipschitz in
  $x$ and $\xi_\theta$ is Lipschitz in $x$ and smooth on the
  $\alpha_\mu$ scale, within the support of $l$ we obtain
\begin{equation}
e \in C^1_\mu S(\alpha^2 \mu,g_\alpha)
\label{amureg}\end{equation}

Combining this with the regularity of the symbol
$\tilde\tphi^{+,\alpha}_{\theta,\mu} \in S(1,g_\alpha) $ we obtain the
symbol regularity for $l$,
\begin{equation}
l \in C^1_\mu S((\alpha^2 \mu)^{-1},g_\alpha)
\label{amurega}\end{equation}

To prove part (a) of the Lemma we observe that for fixed $(t,x)$
the symbol $l(t,x,\xi)$ is a smooth bump function of size
$(\alpha^2 \mu)^{-1}$ in a rectangle of size $ \mu \times (\alpha
\mu)^{n-1}$ oriented in the $ \xi_\theta^{\alpha_\mu}$ direction.
This implies that its kernel $K(t,x,y)$ is bounded by $(\alpha^2
\mu)^{-1}$ times an integrable bump function on the dual scale,
\[
|K(t,x,y)| \lesssim (\alpha^2 \mu)^{-1}  \mu (\alpha \mu)^{n-1}
(1 + \mu |\xi_\theta^{\alpha_\mu}(t,x)(x-y)| +
\alpha \mu |\xi_\theta^{\alpha_\mu}(t,x)\wedge (x-y)|)^{-N}
\]
This bound is symmetric; indeed, since $\xi_\theta^{\alpha_\mu}(t,x)$
is Lipschitz in $x$ we can replace it by $\xi_\theta^{\alpha_\mu}(t,y)$
in the above bound. Thus integrating we have
\[
\sup_x \int |K(t,x,y)| dy \lesssim  (\alpha^2 \mu)^{-1}, \qquad
\sup_y \int |K(t,x,y)| dx \lesssim  (\alpha^2 \mu)^{-1}
\]
 The $L^p$ bounds for $L(t,x,D)$  and also for $\tilde L$ immediately follow.

For later use in the proof we observe that within the support of $l$
we have
\[
|\xi_\theta^{\alpha_\mu}(t,x) \wedge \xi | \lesssim \alpha \mu
\]
Then the same argument as above yields the additional bounds
\begin{equation}
\| (\xi_\theta^{\alpha_\mu}(t,x) \wedge D)^\beta \tilde L u\|_{L^p}
\lesssim  (\alpha \mu)^{|\beta|}  (\alpha^2 \mu)^{-1} \| u\|_{L^p}
\label{wedga}\end{equation}

 For part (b) we write
\[
(A(t,x,D)- \tA(t,x,D)) \tilde L - \tilde\tphi(t,x,D) = R_1(t,x,D) +
R_2(t,x,D)
\]
where
\[
 R_1(t,x,D)=E(t,x,D)
 \tilde S_\mu(D) L(t,x,D) - \tilde \tphi(t,x,D),
\]
respectively
\[
 R_2(t,x,D)=(A_{>\mu^\frac12}(t,x,D)- \tA_{>\mu^\frac12}(t,x,D))
 \tilde S_\mu(D) L(t,x,D),
\]

The operator $R_1$ is localized at frequency $\mu$. The principal part
cancels, and since $e \in C^1_\mu S(\alpha^2 \mu,g_\alpha)$
and $l \in C^1_\mu S((\alpha^2 \mu)^{-1},g_\alpha)$
by the pseudodifferential calculus it follows that
\[
 R_1(t,x,D) \in C^0_\mu S((\alpha^2 \mu)^{-1},g_\alpha)
\]
In addition, the symbol of $R_1$ decays rapidly away from the support
of $\tilde\tphi^{+,\alpha}_{\theta,\mu}$.  Hence we obtain the same kernel
and $L^p$ bounds as in the case of $L(t,x,D)$.

Consider now the operator $R_2$. Since $a(t,x,\zeta)$ is Lipschitz
in $x$ it follows that $|a_{>\mu^\frac12}(t,x,\zeta)| \lesssim
\mu^{-\frac12} |\zeta|$.  Expanding $a_{>\mu^\frac12}(t,x,\zeta)$
in a rapidly decreasing series of spherical harmonics with respect
to $\zeta$, we can separate variables and reduce the problem to
the simpler case when $a_{>\mu^\frac12}(t,x,\zeta) = b(t,x)
c(\zeta)$ with $|b| < \mu^{-\frac12}$ and $c$ is smooth and
homogeneous of order $1$. For the symbol $c-\tc$ we use the
representation
\[
c(\zeta) - \tc(t,x,\zeta) = \psi(\xi_\theta^{\alpha_\mu}, \zeta)
(\xi_\theta^{\alpha_\mu}(t,x) \wedge \zeta)^2
\]
where $\psi$ is smooth in both arguments and homogeneous of order
$-1$ in $\zeta$. Separating variables in $\psi$ we can assume without
any restriction in generality that $\psi$ depends only on $\zeta$.
Then after some simple commutations we obtain
\[
c(D) - \tc(t,x,D) = (\xi_\theta^{\alpha_\mu}(t,x) \wedge D)^2
\psi(D) + O(1)_{L^p \to L^p}
\]
To estimate this we use \eqref{wedga}.
The factor $\psi(D)\tilde S_\mu(D)$ yields
an extra $\mu^{-1}$ factor in the $L^p$ bounds, therefore we obtain
\[
\| R_2(t,x,D)\|_{L^p \to
  L^p} \lesssim \mu^{-\frac12}
\]


Finally we prove part (c). By (a), $\tilde L$ is $L^2$ bounded with norm
$O(\alpha^{-2} \mu^{-1})$, therefore
it remains to prove the commutator estimate
\begin{equation}
\| [ D_t+ A_{<\mu^\frac12}(t,x,D),\tilde S_\mu  L(t,x,D)]\|_{L^\infty L^2 \to L^2}
\lesssim \alpha^{-2} \mu^{-1}
\label{lcom}\end{equation}
This is a consequence of the operator bound
\[
 [ D_t+ A_{<\mu^\frac12}(t,x,D),\tilde S_\mu  L(t,x,D)] \in L^2 C^0_\mu S(\alpha^{-2} \mu^{-1},g_\alpha)
\]
To prove it we use the pdo calculus to represent the commutator
as a principal term plus a second order error,
\[
[ D_t+ A_{<\mu^\frac12}(t,x,D),\tilde S_\mu  L(t,x,D)] = \tilde S_\mu Q(t,x,D) + R(t,x,D)
\]
where the  principal part $q$ has symbol
\[
q(t,x,\xi) =-i \{\tau + a_{<\mu^\frac12}(t,x,\xi), l(t,x,\xi)\}
\]
The remainder $R$ is localized at frequency $\mu$.
A direct computation, using \eqref{amurega}, shows that  its symbol
satisfies
\[
r \in L^2 C^0_\mu S(\alpha^{-2} \mu^{-1},g_\alpha)
\]
It remains to consider the above Poisson bracket and prove that
\begin{equation}
q \in L^2 C^0_\mu S(\alpha^{-2} \mu^{-1},g_\alpha)
\label{qsal}\end{equation}
For this we write $q$ in the form
\[
iq  = - \tilde\tphi^{+,\alpha}_{\theta,\mu} q_1 e^{-2}   + q_2 e^{-1}
+ q_3 e^{-1}
\]
where
\[
q_1 (t,x,\xi)= \left\{\tau +
  a_{<\mu^\frac12} , e \right\}, \qquad
q_2 (t,x,\xi)= \left\{\tau + a_{<\alpha^{-1}}, \tilde\tphi^{+,\alpha}_{\theta,\mu}\right\}
\]
respectively
\[
q_3 (t,x,\xi)=
\left\{ a_{\alpha^{-1}<\cdot < \mu^{\frac12}}, \tilde\tphi^{+,\alpha}_{\theta,\mu}\right\}
\]

Within the support of $\tilde\tphi^{+,\alpha}_{\theta,\mu}$ we
know that $e \in C^1_\mu S(\alpha^2 \mu, g_\alpha)$ is an elliptic
symbol. Hence for the first term it suffices to show that $q_1 \in
C^0_\mu S(\alpha^2 \mu, g_\alpha)$.  Indeed, by definition $q_1$
is a homogeneous symbol of order $1$ which is continuous in $x$
and homogeneous in $\zeta$. In addition, we know that
$e(t,x,\zeta)$ vanishes of second order in $\zeta$ at $
(t,x,\xi_\theta^{\alpha_\mu}(x,t))$ which is also invariant with
respect to the $\tau+ a_{<\mu^\frac12}$ Hamilton flow. Then $q$
must vanish of second order in $\zeta$ at
$(t,x,\xi_\theta^{\alpha_\mu}(x,t))$. Arguing as in the case of
$e$, this implies that within the support of
$\tilde\tphi^{+,\alpha}_{\theta,\mu}$ we have $q_1 \in C^0_\mu
S(\alpha^2 \mu, g_\alpha)$.

As in \eqref{lowcoma} we know that $q_2 \in S(1,g_\alpha)$.
Also we have $a_{\alpha^{-1}<\cdot < \mu^{\frac12}} \in L^2 C^2_\mu
S(\alpha^2 \mu, g_\alpha)$ and  $\tilde\tphi^{+,\alpha}_{\theta,\mu}
\in S(1,g_\alpha)$ therefore $q_3 \in L^2 C^1_\mu S(1,g_\alpha)$.

This concludes the proof of \eqref{qsal} and therefore the proof of
the lemma.

\end{proof}

To continue the estimate of term II in Case 1 we define the auxiliary
trilinear form
\begin{eqnarray*}
E(u,v,\tw) &=& \int (D_t + A(t,x,D)) u\, v \tw dx dt +  \int  u (D_t -
A(t,x,-D)) v\, \tw dx dt
\\ &+&  \int  u v\, (D_t + \tA(t,x,D)) \tw dx dt
\end{eqnarray*}
With $\tilde w= \tilde L  w_{\theta}^\alpha$ we
write
\begin{equation}
\begin{split}
\int u_\theta v_\theta w_\theta^\alpha dx dt =& -\int u_\theta
v_\theta\, ((A(t,x,D)- \tA(t,x,D))\tilde L-1)w_\theta^\alpha dx dt
\\ +& \int (D_t + A(t,x,D)) u_\theta\, v_\theta \tw dx dt  \\ +& \int  u_\theta (D_t - A(t,x,-D)) v_\theta\,
\tw dx dt
\\ +& \int  u_\theta v_\theta (D_t + A(t,x,D)) \tw dx dt
\\ -& E(u_\theta,v_\theta,\tw)
\end{split}
\label{longsum}
\end{equation}
We bound each term separately. For the first one we write
\[
\begin{split}
(A(t,x,D)- \tA(t,x,D)) \tilde{ L}-1 =& (A(t,x,D)-\tA(t,x,D))\tilde
L-\tilde \tphi^{+,\alpha}_{\theta,\mu}(t,x,D)
\\ & + (\tilde \tphi^{+,\alpha}_{\theta,\mu}(t,x,D)-1)
\end{split}
\]
The contribution of the first line is estimated using Lemma~\ref{L}
(b) and \eqref{pecher} for $w_\theta^\alpha$,
\[
\begin{split}
\bigg| \int  u_\theta v_\theta  [&(A(t,x,D)-\tA(t,x,D))  \tilde L
- \tilde \tphi^{+,\alpha}_{\theta,\mu}(t,x,D)]w_\theta^\alpha dx
dt \bigg|
\\& \lesssim (\alpha^{-2}\mu^{-1}+\mu^{-\frac12})  \|u_\theta\|_{L^\infty L^2} \|v_\theta\|_{L^\infty L^2} \|w_\theta^\alpha \|_{L^{2}
L^{\infty}}
\\& \lesssim (\alpha^{-2}\mu^{-1} +\mu^{-\frac12}) \|u_\theta\|_{X_+} \|v_\theta\|_{X_-} \|w_\theta^\alpha \|_{L^{2}
L^{\infty}}
\\& \lesssim (\alpha^{-2}\mu^{-1}+\mu^{-\frac12}) \alpha^\frac12 \mu^\frac32
\|u_\theta\|_{X_+} \|v_\theta\|_{X_-}
\|w_\theta^\alpha\|_{X_+^{\mu,\alpha,\theta}}
\end{split}
\]
For the contribution of the second line we observe that
\[
(\tilde  \tphi^{+,\alpha}_{\theta,\mu}(t,x,D)-1) w_\theta^\alpha = (\tilde \tphi^{+,\alpha}_{\theta,\mu}(t,x,D)-1) \tilde
\phi_\theta^{+,\alpha}  S_\mu w
\]
where the symbols $\tilde \tphi_{\theta,\mu}^{+,\alpha} -1$ and
$\tilde \phi_{\theta,\mu}^{+,\alpha} s_\mu $ have disjoint supports.
Since they both belong to $S(1,g_\alpha)$, this yields a gain of a
factor $(\alpha^2 \mu)^{-N}$ in \eqref{pecher}, with $N$ arbitrarily
large:
\[
\sum_\theta \| (\tilde \tphi^{+,\alpha}_{\theta,\mu}(t,x,D) -1) w_\theta^\alpha\|_{L^2 L^\infty}^2 \lesssim
\mu^{\frac52} (\alpha^2 \mu)^{-N} \|w_\theta^\alpha\|_{X_+^{\mu,\alpha,\theta}}^2
\]
This is more than we need.

For the second term in \eqref{longsum} we use the $L^2$ bound for
$(D_t+A)u_\theta$, the energy bound for $v_\theta$ and
\eqref{pecher} for $\tw$. This yields
\[
\left| \int (D_t +A(t,x,D)) u_\theta\, v_\theta \tw dx dt \right|
\lesssim
 \alpha^{-\frac32} \mu^{\frac12} \|u_\theta\|_{X_+} \|v_\theta\|_{X_-}
\|w_\theta^\alpha\|_{X_+^{\mu,\alpha,\theta}}
\]
The third term is similar.

For the fourth term in  \eqref{longsum} we use the energy for the
first two factors combined with Bernstein derived $L^2 L^\infty$ bound
for the third,
\[
\begin{split}
\left|  \int  u_\theta v_\theta\, (D_t + A(t,x,D)) \tw dx dt\right| &
\lesssim \|u\|_{X_+} \|v\|_{X_-}
 \|(D_t + A(t,x,D)) \tw \|_{L^2 L^\infty} \\ &\lesssim (\alpha^2
 \mu)^{-1} (\mu (\alpha
 \mu)^3)^\frac12 \|u_\theta\|_{X_+} \|v_\theta\|_{X_-}
\|w_\theta^\alpha\|_{X_+^{\mu,\alpha,\theta}}
\end{split}
\]

It remains to prove  the estimate  for $E$. Observe that the time
derivatives in $E$ can be integrated out, producing
contributions of the form
\begin{equation}
\int u_\theta\, v_\theta \, \tw dx
\label{puvtw}\end{equation}
at the initial and the final time.  These are estimated using energy bounds
for the first two factors and the pointwise bound arising from
Bernstein's inequality for the last factor,
\[
\| \tw\|_{L^\infty} \lesssim (\alpha^2 \mu)^{-1} \|
w^\alpha_\theta\|_{L^\infty} \lesssim (\alpha^2 \mu)^{-1} (\mu
(\alpha \mu)^3)^\frac12
\|w^\alpha_\theta\|_{X_+^{\mu,\alpha,\theta}} = (\alpha^2
\mu)^{-\frac14}
 \mu^\frac54\|w^\alpha_\theta\|_{X_+^{\mu,\alpha,\theta}}
\]

This leaves us with a purely spatial trilinear form,
\[
\int  E_0(u_\theta,v_\theta,\tw) dt
\]
where
\[
E_0(u,v,\tw) = \int A(t,x,D) u\, v \tw  - u A(t,x,-D) v \,\tw +  u
v \, \tA(t,x,D) \tw \ dx
\]
 The main bound for $E_0$ is provided in the next lemma.

\begin{lemma} \label{leob}
  Let $1 \leq \mu \lesssim \lambda$. Assume that $\xi_\theta$ is a
  Lipschitz function of $x$ with $|\xi_\theta -\theta| \ll 1$ and that $a \in C^1 S^1_{hom}$.  Then the
  trilinear form $E_0$ satisfies the fixed time estimate:
\begin{equation}
\begin{split}
|E_0(u,v,\tw)| \lesssim& \  \|u\|_{L^{p_1}}
\|v\|_{L^{q_1}}  \|\tw\|_{L^{r_1}} \\ &  + \lambda^{-1} \|   (\xi_\theta \wedge D) u\|_{L^{p_2}}
\|v\|_{L^{q_2}}  \| (\xi_\theta \wedge D)  \tw\|_{L^{r_2}}
\\ & +   \lambda^{-1} \|  u\|_{L^{p_2}}
\| (\xi_\theta \wedge D) v\|_{L^{q_2}}  \| (\xi_\theta \wedge D)  \tw\|_{L^{r_2}}
\\ & + \mu \lambda^{-2} \|  (\xi_\theta \wedge D) u\|_{L^{p_3}}
\| (\xi_\theta \wedge D) v\|_{L^{q_3}}  \|\tw\|_{L^{r_3}}
\end{split}
\label{eob}\end{equation}
for all indices
\[
\frac{1}{p_i} + \frac{1}{q_i} + \frac{1}{r_i} = 1, \qquad 1 \leq
p_i,q_i,r_i \leq \infty
\]
and for all functions $u$, $v$ localized at frequency $\lambda$ in a small
angular neighbourhood of $\theta$, respectively $-\theta$ and all $w$
localized at frequency $\mu$.
\label{ee} \end{lemma}
While any choice of $L^p$ norms is allowed in the lemma, in order to
conclude the proof of the estimate for $E$ it suffices to use the set
of indices $(2,2,\infty)$. We apply the lemma with $u = u_\theta$,
$v = v_\theta$ and $\tw = \tilde L w^\alpha_\theta$ as above.
This yields
\begin{equation}
\begin{split}
\left|\int  E_0(u_\theta,v_\theta,\tw) dt\right| \lesssim & \
 \|u_\theta\|_{L^{\infty} L^2}
\|v_\theta\|_{L^{\infty} L^2}  \|\tw\|_{L^{2} L^\infty}  \\ & + \lambda^{-1} \|   (\xi_\theta \wedge D) u_\theta\|_{L^{\infty} L^2}
\|v_\theta\|_{L^{\infty} L^2}  \| (\xi_\theta \wedge D)  \tw\|_{L^{2} L^\infty}
\\ & +   \lambda^{-1} \|  u_\theta\|_{L^{\infty} L^2}
\| (\xi_\theta \wedge D) v_\theta\|_{L^{\infty} L^2}  \| (\xi_\theta \wedge D)  \tw\|_{L^{2} L^\infty}
\\ & + \mu \lambda^{-2} \|  (\xi_\theta \wedge D) u_\theta\|_{L^{\infty} L^2}
\| (\xi_\theta \wedge D) v_\theta\|_{L^{\infty} L^2}  \|\tw\|_{L^{2} L^\infty}
\end{split}
\label{secondee}\end{equation}
Due to the angular localization, the operator $ (\xi_\theta \wedge D) $
yields a factor of $\mu^{-\frac12} \lambda$ when applied to $u_\theta$ or
$v_\theta$, respectively a factor of $\alpha \mu$ when applied to
$\tilde w$.  Hence we obtain
\[
\left|\int E_0(u_\theta,v_\theta,\tw) dt\right| \lesssim
\frac{\mu^{\frac32} \alpha^\frac12}{\alpha^2 \mu}  (1+ \alpha
\mu^\frac12 +\alpha \mu^\frac12+ 1)
\|u_\theta\|_{X_+^{\l,\alpha,\theta}}
\|v_\theta\|_{X_-^{\l,\alpha,\theta}}
\|w^\alpha_\theta\|_{X_+^{\mu,\alpha,\theta}}
\]
which is acceptable since $\alpha^2 \mu \geq 1$.

\begin{proof}[Proof of Lemma~\ref{leob}:]
  Since the symbol $a$ is smooth and homogeneous of order $1$ with
  respect to $\xi$, we can use its representation in terms of the
  spherical harmonics and reduce the problem to the case when $a$ has
  the form
\[
a(x,\xi) = b(x) c(\xi)
\]
where $b$ is Lipschitz continuous.

We denote by $\xi$, respectively $\eta$ the frequencies for the
$u_\theta$, respectively $v_\theta$ factors in $E_0$. Then $\xi$ and $\eta$
have size $\lambda$ and are in a small angular neighbourhood of $\theta$.
 We expand $c$ around the line generated by
$\xi_\theta$ into a linear term and a quadratic error,
\[
c(\xi) = \xi  (\nabla c)(\xi_\theta) + \xi B(\xi,\xi_\theta) \xi
\]
where $B$ is homogeneous of  order $-1$ with respect to $\xi$ and
can be chosen so that
\[
\xi_\theta B(\xi,\xi_\theta) = 0, \qquad  B(\xi,\xi_\theta) \xi_\theta = 0
\]
To see that this is possible we observe that after a rigid rotation we
can assume that $\xi_\theta = e_1$. For $\xi = (1,\xi')$ with $|\xi'|
\ll 1$ we write the first order Taylor polynomial with integral
remainder
\[
\begin{split}
c(1,\xi') = &\  c(1,0) + \xi' c_{\xi'}(1,0) + \xi'  B(1,\xi') \xi'
\\ = &\  c_{\xi_1}(1,0) + \xi' c_{\xi'}(1,0) + \xi'  B(1,\xi') \xi'
\end{split}
\]
where $B$ is given by
\[
 B(1,\xi') = \int_0^1 (1-h) \nabla^2_{\xi'} a(1,h\xi') dh
\]
This extends by homogeneity to all $\xi$ in a small angle around
$\theta$.

We represent $B$ as a rapidly convergent sum of terms of the form
\[
\lambda^{-1} F(\xi_\theta) g(\xi)
\]
where $g$ is a scalar function which is bounded and smooth on the
$\lambda$ scale and $F$ is a matrix inheriting the above property of
$B$,
\begin{equation}
\xi_\theta F(\xi_\theta) = 0, \qquad F(\xi_\theta) \xi_\theta = 0
\label{null}\end{equation}
 So we have
\[
c(\xi) = \xi  (\nabla c)(\xi_\theta) + \lambda^{-1} \sum
 \xi  F(\xi_\theta) \xi g(\xi)
\]
Then we obtain the rapidly convergent series representation
\begin{eqnarray*}
c(\xi) - c(\eta) &=& (\xi-\eta) (\nabla c)(\xi_\theta) +\lambda^{-1}  \sum (\xi - \eta)
F(\xi_\theta) \xi g(\xi) \\ &+& \lambda^{-1} \sum \eta F(\xi_\theta) (\xi -\eta) g(\eta)
\\ &+& \lambda^{-2} \sum \eta F(\xi_\theta) \xi  (\xi - \eta)  h(\xi) k(\eta)
\end{eqnarray*}
where $h$ and $k$ are smooth and bounded on the $\lambda$ dyadic
scale.

We use this representation for the first two components in $E_0$.  The
contribution of the first term above cancels the principal part of the
third component in $E_0$.  We retain the other three terms though,
therefore this yields the following rapidly convergent series
representation for $E_0$:
\[
E_0(u,v,\tw) =    \int    u v  \tw D(b (\nabla c)(\xi_\theta)) dx
+ \sum E_0^1 + \sum E_0^2 + \sum E_0^3
\]
The first term is easily estimated since $b (\nabla
c)(\xi_\theta)$ is Lipschitz continuous. The first summand has the
form
\[
\begin{split}
E_0^1 &= \l^{-1} \int   F(\xi_\theta) D (Dg(D) u v)\, \tw dx
\\ & =  - \l^{-1} \int  D F(\xi_\theta) D g(D) u v  \tw
+ F(\xi_\theta) D g(D) u v  D \tw \ dx
\end{split}
\]
In the first term $F(\xi_\theta)$ is Lipschitz in $x$ and the $u$
derivative yields a factor of $\lambda$.  For the second term on the
other hand we use \eqref{null} to estimate
\[
| D g(D) u F(\xi_\theta) D \tw| \lesssim   |(\xi_\theta \wedge D) g(D) u|  |(\xi_\theta \wedge D) \tw|
\]
 Commuting $g(D)$ with $(\xi_\theta \wedge D)$ we get
\[
|(\xi_\theta \wedge D) g(D) u| \leq |g(D) (\xi_\theta \wedge D) u| + | [
g(D), \xi_\theta \wedge D] u|
\]
with the commutator $[ g(D), \xi_\theta \wedge D]$ bounded in all
$L^p$ spaces. Hence
\[
|E_0^1| \lesssim \| u\|_{L^{p_1}} \|v\|_{L^{q_1}} \|\tw\|_{L^{r_1}} +
\lambda^{-1} \| (\xi_\theta \wedge D) u\|_{L^{p_2}} \|v\|_{L^{q_2}}
\|(\xi_\theta \wedge D) \tw\|_{L^{r_2}}
\]
The second summand of $E_0$ is similar but with the roles of $u$ and
$v$ reversed.

Finally,
\[
\begin{split}
E_0^3 &=  \lambda^{-2} \int F(\xi_\theta)\, D (D h(D) u \,D k(D) v)\, \tw dx
\\ &=  -\lambda^{-2} \int D F(\xi_\theta) \,D h(D) u\, D k(D) v \, \tw
+ F(\xi_\theta)\, D h(D) u\, D k(D) v \, D \tw \ dx
\end{split}
\]
where the matrix $F(\xi_\theta)$ is paired with the $u$ and $v$ derivatives.
In the first term the  two derivatives on $u$ and $v$ yield a
$\lambda^2$ factor. In the second term we use  as before \eqref{null}
and commute out the $h(D)$ and $k(D)$ multipliers.
We obtain
\[
|E_0^3| \lesssim\| u\|_{L^{p_1}} \|v\|_{L^{q_1}} \|\tw\|_{L^{r_1}} + \mu
\lambda^{-2}\| (\xi_\theta \wedge D) u\|_{L^{p_3}} \| (\xi_\theta \wedge D) v\|_{L^{q_3}} \|
\tw\|_{L^{r_3}}
\]
Summing up the results we get the conclusion of the Lemma.

\end{proof}

{\bf Case 1, term III}:
This has the form
\[
III = \int  \sum_{ \alpha > \mu^{-\frac12} }\sum_{\theta \in
  O_{\alpha}}
\phi^{+,\alpha}_{\theta,\lambda}(t,x,D) S_\l u \ \tilde
\phi^{-,\alpha}_{-\theta,\lambda}(t,x,D) S_\l v \ \tilde
\phi^{+,\alpha \mu^{-1} \lambda}_{\theta,\mu}(t,x,D) S_\mu w \,dx
dt
\]
In this case the summation with respect to  $\theta$ is
accomplished by \eqref{limsumeqa}, while for the $\alpha$ summation
we simply accept a $\ln \mu$ loss.  Fixing $\alpha$ and $\theta$ we set
\[
u_\theta^\alpha = \phi^{+,\alpha}_{\theta,\lambda}(t,x,D) S_\l u, \quad v_\theta^\alpha =
\phi^{-,\alpha}_{-\theta,\lambda}(t,x,D) S_\l v, \quad w_\theta^\alpha=\tilde \phi^{+,\alpha
\mu^{-1} \lambda}_{\theta,\mu}(t,x,D) S_\mu w.
\]
and repeat the analysis for Case 1, term II. The angular localization
of $u_\theta^\alpha $ and $v_\theta^\alpha$ is not used in the bounds
for the first four terms in \eqref{longsum}, therefore that part of
the argument rests unchanged. The same applies to the bound for the
fixed time integral in \eqref{puvtw}.

 It remains to consider the bound for $E(u_\theta^\alpha,v_\theta^\alpha,\tw)$.
The $\alpha$ localization angle for $w_\theta^\alpha$ is now $\alpha
\mu^{-1} \lambda$, therefore part (b) of Lemma~\ref{L} gives
\[
\| \tw\|_{X_+} \lesssim \frac{\mu}{\alpha^2 \lambda^2}
\|w_\theta^\alpha\|_{X_+}
\]
This is stronger than in the previous case because it gives a high
frequency gain. Now we are able to use Lemma~\ref{ee} with exponents
$(3,2,6)$  to obtain
\[
\begin{split}
  |\int E_0(u_\theta^\alpha,v_\theta^\alpha,\tw) dt| \lesssim &\
  \|u_\theta^\alpha\|_{L^2 L^3} \|v_\theta^\alpha\|_{L^{\infty} L^2}
  \|\tw\|_{L^{2} L^6}
\\ & + \lambda^{-1} \| (\xi_\theta \wedge D)  u_\theta^\alpha\|_{L^2 L^3} \|v_\theta^\alpha\|_{L^{\infty} L^2}
\| (\xi_\theta  \wedge D) \tw\|_{L^{2} L^6} \\ & + \lambda^{-1} \|
  u_\theta^\alpha\|_{L^2 L^3} \| (\xi_\theta \wedge D)
  v_\theta^\alpha\|_{L^{\infty} L^2} \| (\xi_\theta \wedge D) \tw\|_{L^{2}
    L^6} \\ & + \mu \lambda^{-2} \| (\xi_\theta \wedge D)
  u_\theta^\alpha\|_{L^2 L^3} \| (\xi_\theta \wedge D)
  v_\theta^\alpha\|_{L^{\infty} L^2} \|\tw\|_{L^{2} L^6}
\end{split}
\]
Due to the angular localization on the $\alpha$ scale for
$u_\theta^\alpha$ and $v_\theta^\alpha$, respectively on the
$\alpha \mu^{-1} \lambda$ scale for $w_\theta^\alpha$, all $(\xi_\theta \wedge D)$
operators above yield $\alpha \lambda$ factors. Hence, taking advantage
of the Strichartz estimates, we obtain
\[
\begin{split}
 |\int E_0(u_\theta^\alpha,v_\theta^\alpha,\tw) dt| \lesssim&\
 \frac{\mu}{\alpha^2 \lambda^2} \alpha^2 \lambda\  \lambda^{\frac5{12}} \mu^\frac56
 \|u_\theta^\alpha\|_{X_+^{\l,\alpha,\theta}} \|v_\theta^\alpha\|_{X_-^{\l,\alpha,\theta}}
\|w_\theta^\alpha\|_{X_+^{\mu,\frac{\alpha\l}{\mu},\theta}}
\\
=&\ \lambda^{-\frac7{12}} \mu^\frac{11}6
 \|u_\theta^\alpha\|_{X_+^{\l,\alpha,\theta}} \|v_\theta^\alpha\|_{X_-^{\l,\alpha,\theta}}
\|w_\theta^\alpha\|_{X_+^{\mu,\frac{\alpha\l}{\mu},\theta}}
\end{split}
\]
which is satisfactory since $\lambda \gtrsim \mu$.

We conclude this case with two remarks. First, in this context the
proof of Lemma~\ref{ee} is somewhat of an overkill. In fact, it would
suffice to linearize separately $a(t,x,\xi)$ and $a(t,x,\eta)$ around $\xi_\theta$
and use the fact that the symbol $a(t,x,\xi) - \ta(t,x,\xi)$ has size
$\alpha^2 \lambda$ at frequency $\l$ in $H_\alpha S_\alpha(\theta)$.
Secondly, the endpoint Strichartz estimate is only used here for
convenience; there is some flexibility in choosing the indices.

{\bf Case 1, term IV.}
This has the form
\[
IV = \int  \sum_{ \alpha > \mu^{-\frac12} }\sum_{\theta \in
  O_{\alpha}}
\phi^{+,\alpha}_{\theta,\lambda}(t,x,D) S_\l u \ \tilde
\phi^{-,\alpha}_{-\theta,\lambda}(t,x,D) S_\l v \ \phi^{+,\alpha
\mu^{-1} \lambda}_{\theta,\mu}(t,x,D) S_\mu w \,dx dt
\]
Again the summation with respect to $\theta$ is accomplished by
\eqref{limsumeqa}, while for the $\alpha$ summation we simply accept a
$\ln \mu$ loss. This term  is better behaved because the symbol
\[
\phi^{+,\alpha}_{\theta,\lambda}(x,\xi)\, \tilde
\phi^{-,\alpha}_{-\theta,\lambda}(x,\eta) \,\phi^{+,\alpha \mu^{-1}
\lambda}_{\theta,\mu}(x,\zeta)
\]
vanishes on $H=\{\xi+\eta + \zeta = 0\}$.  Precisely,
in the support of the above symbol we have
\[
|\xi| \approx \lambda,\ |\xi \wedge \xi_\theta^\alpha| \lesssim \alpha
\lambda, \qquad |\eta| \approx \lambda,\ |\eta \wedge
\xi_\theta^\alpha| \approx C \alpha
\lambda, \qquad  |\zeta| \approx \lambda,\ |\zeta \wedge
\xi_\theta^\alpha| \lesssim \alpha
\lambda.
\]
This leads to
\begin{equation}
|(\xi + \eta +\zeta) \wedge
\xi_\theta^\alpha| \approx C \alpha
\lambda
\label{offdiag}\end{equation}
This can be taken advantage of in a direct computation in the above
formula. Including the dyadic frequency localizations into the
$\phi$'s, each term in $IV$ has the integral representation
\[
\int \phi^{+,\alpha}_{\theta,\lambda}(t,x,\xi)  \hat u(\xi) \
\tilde \phi^{-,\alpha}_{-\theta,\lambda}(t,x,\eta)  \hat v(\eta) \
\phi^{+,\alpha \mu^{-1} \lambda}_{\theta,\mu}(t,x,\zeta) \hat
w(\zeta)\, e^{i x(\xi+\eta+\zeta)} \, d\xi d\eta d\zeta dx dt
\]
Defining the spatial elliptic operator $F$ with symbol
\[
f(t,x,\xi) =  (\xi \wedge \xi_\theta^\alpha)^{2N}
\]
we have
\[
F(t,x,D_x)  e^{i x(\xi+\eta+\zeta)} = |(\xi + \eta +\zeta) \wedge
\xi_\theta^\alpha|^{2N} e^{i x(\xi+\eta+\zeta)}
\]
Hence integration by parts in the above formula leads to
\[
\int \psi(t,x,\xi,\eta,\zeta)   \hat u(\xi) \, \hat v(\eta)  \, w(\zeta) e^{i
x(\xi+\eta+\zeta)} \, d\xi d\eta d\zeta dx dt
\]
where the new symbol $\psi$ is
\[
\psi(t,x,\xi,\eta,\zeta) = F^*(t,x,D_x) \left( \frac{\phi^{+,\alpha}_{\theta,\lambda}(t,x,\xi) \tilde
\phi^{-,\alpha}_{-\theta,\lambda}(t,x,\eta)   \phi^{+,\alpha
\mu^{-1} \lambda}_{\theta,\mu}(t,x,\zeta) }{  |(\xi + \eta +\zeta) \wedge
\xi_\theta^\alpha|^{2N}} \right)
\]
In the support of the numerator the bound \eqref{offdiag}
holds. Hence separating the variables we can represent the denominator
as a rapidly convergent series with terms
\[
(\alpha \lambda)^{-2N} \chi_{< \alpha \l }( \xi \wedge \xi_\theta^\alpha) \chi_{C \alpha \l }
( \eta \wedge \xi_\theta^\alpha)
\chi_{< \alpha \l }( \zeta \wedge \xi_\theta^\alpha)
\]
where each of the $\chi$'s above is a unit bump function on the
$\alpha \lambda$ scale. Thus they can be included in the corresponding
$\phi$ factors. Due to the $S(g_\alpha)$ regularity of the $\phi$
factors, each derivative $ \xi_\theta^\alpha \wedge D$ applied to them
yields an $\alpha^{-1}$ factor. Thus $\psi$ is represented
as a rapidly convergent series of products of the form
\[
(\alpha^2 \lambda)^{-2N}  \psi^{+,\alpha}_{\theta,\lambda}(t,x,\xi) \tilde
\psi^{-,\alpha}_{-\theta,\lambda}(t,x,\eta)   \psi^{+,\alpha
\mu^{-1} \lambda}_{\theta,\mu}(t,x,\zeta)
\]
where the $\psi$ factors have the same support and regularity as the
corresponding $\phi$'s.  The integral above is similarly represented
as a rapidly convergent series with terms of the form
\[
(\alpha^2 \lambda)^{-2N} \int \psi^{+,\alpha}_{\theta,\lambda}(t,x,D) S_\lambda u \,\tilde
\psi^{-,\alpha}_{-\theta,\lambda}(t,x,D) S_\l v \, \psi^{+,\alpha
\mu^{-1} \lambda}_{\theta,\mu}(t,x,D) S_\mu w \,dx dt
\]
Since $\alpha > \mu^{-\frac12}$, the factor in front of the above
integral allows us to exchange low frequencies for high frequencies.
This suffices in order to bound the last integral using Strichartz
estimates.

{\bf Case 1, term V}

This is similar to Case 1, term $IV$. This time in the support of
the symbol
\[
\phi^{+,\alpha}_{\theta,\l}(t,x,\xi)\, \tilde
\phi^{-,\alpha}_{-\theta,\l}(t,x,\eta)\,
\tilde\phi^{+,\beta}_{\theta,\mu}(t,x,\zeta)
\]
we have
\[
|\xi| \approx \lambda,\ |\xi \wedge \xi_\theta^\alpha| \lesssim \alpha
\lambda, \qquad |\eta| \approx \lambda,\ |\eta \wedge
\xi_\theta^\alpha| \approx C \alpha
\lambda, \qquad  |\zeta| \approx \lambda,\ |\zeta \wedge
\xi_\theta^\alpha| \approx C\beta
\lambda.
\]
Hence
\[
|(\xi + \eta +\zeta) \wedge
\xi_\theta^\alpha| \approx C \alpha
\lambda
\]
therefore the symbol above is supported at distance $\beta \l$ from
the diagonal $H$.  Hence integrating by parts as in the previous case
we gain arbitrary powers of $(\alpha \beta \l)^{-1}$. Then we can
close the argument using Strichartz type estimates.

{\bf Case 2}, $1 < d < \mu$. This requires only minor changes,
which we describe in what follows.  We still consider the five
terms in the trilinear decomposition \eqref{trilindec}, but we
replace the smallest localization angle $\mu^{-\frac12}$ by
$d^{\frac12} \mu^{-\frac12}$.

{\bf Case 2, term I}. Here we need \eqref{limsumeqa} to sum
expressions of the form
\[
I =  \int \phi^{+,d^{\frac12}\mu^{-\frac12}}_{\theta,\lambda}(t,x,D) S_\l u \
\phi^{-,d^{\frac12}\mu^{-\frac12}}_{-\theta,\lambda}(t,x,D) S_\l v\
 \phi^{+,d^{\frac12}\mu^{-\frac12}}_{\theta,\mu}(t,x,D) S_\mu w \,dx
\]
over $\theta \in O_{d^{\frac12}\mu^{-\frac12}}$.  Each term is bounded
by combining the energy estimate for the first factor, the $L^4 L^2$
bound for the second and \eqref{pecher} for the third.

{\bf Case 2, term II.} Here we use \eqref{limsumeqa}
for the summation of expressions of the form
\[
II =  \int \phi^{+,d^{\frac12}\mu^{-\frac12}}_{\theta,\lambda}(t,x,D) S_\l u\
 \phi^{-,d^{\frac12}\mu^{-\frac12}}_{-\theta,\lambda}(t,x,D) S_\l v\,
\sum_{\alpha > d^{\frac12}\mu^{-\frac12}} \tilde
\phi^{+,\alpha}_{\theta,\lambda}(t,x,D) S_\mu w dx
\]
over $\theta \in O_{d^{\frac12}\mu^{-\frac12}}$. We use the same
operator $L$, the same function $\tilde w$ and the same trilinear form
$E$. In \eqref{longsum} the first, second and fourth terms are
estimated in the same way, but using the $L^4 L^2$ bound for the
second factor. In the third term we lose a power of $d$,
\begin{eqnarray*}
\left|\int u_\theta (D_t - A(t,x,-D)) v_\theta \tw dx\right| &\lesssim&
\|u_\theta\|_{L^\infty L^2} \|(D_t - A(t,x,-D)) v_\theta\|_{L^2}
 \|\tw\|_{L^2 L^\infty}\\ & \lesssim&
\|u_\theta\|_{X_+} d^{\frac34} \|v_\theta\|_{X_{-,d}} \frac{1}{\alpha^2 \mu}
\alpha^{\frac12} \mu^{\frac32} \|w_\theta^\alpha\|_{X_+^{\mu,\alpha,\theta}}
\\ &\lesssim & \left(\frac{d}{\alpha^2 \mu}\right)^\frac34  \mu^{\frac54}
\|u_\theta\|_{X_+}\|v_\theta\|_{X_{-,d}} \|w_\theta\|_{X_+^{\mu,\alpha,\theta}}
\end{eqnarray*}
But this is still acceptable due to the reduced range for $\alpha$,
namely $\alpha^2 \mu \geq d$.

In the expression \eqref{puvtw} there is a $d^{\frac14}$ loss in the
$L^2$ bound for $v_\theta$, but this is compensated for by the
previously unused $(\alpha^2 \mu)^{-\frac14}$ factor in the pointwise
bound for $\tw$.

Finally, for the $E_0$ bounds we reuse \eqref{secondee} but with all
the $v_\theta$ factors estimated in $L^2$. This produces an extra
$d^{-\frac14}$ gain. On the other hand, the angular localization for
$u_\theta$ and $v_\theta$ is worse. Precisely, the operator $
(\xi_\theta \wedge D) $ yields a factor of $d^\frac12 \mu^{-\frac12}
\lambda$ when applied to $u_\theta$ or $v_\theta$, respectively a
factor of $\alpha \mu$ when applied to $\tilde w$.  Hence we obtain
\[
\left|\int E_0(u_\theta,v_\theta,\tw) dt\right| \lesssim
\frac{\mu^{\frac32}  \alpha^\frac12}{d^{\frac14} \alpha^2 \mu}  (1+
d^\frac12 \alpha
\mu^\frac12 +d^\frac12 \alpha
\mu^\frac12+ d) \|u_\theta^\alpha\|_{X_+^{\l,\theta,\alpha}}
\|v_\theta^\alpha\|_{X_-^{\lambda,\theta,\alpha}}
\|w_\theta^\alpha\|_{X_+^{\mu,\theta,\alpha}}
\]
This is still acceptable since $\alpha^2 \mu \geq d$.

{\bf Case 2, term III.} Compared to the similar argument in Case 1,
the following modifications are required:

(i) The third term in \eqref{longsum} is treated as in  Case 2, term
II.

(ii) In the bound for $E_0$, the $L^\infty L^2$ norms are replaced by
$L^4 L^2$ in all the $v_\theta$ factors.

{\bf Case 2, terms IV,V.} These are identical to Case 1.

\end{proof}

\textbf{Acknowledgement} \vspace{0.1in}

Both authors would like to thank MSRI for the hospitality in the
Fall 2005 semester, where part of this article was written. Both
authors were supported in part by the NSF grant DMS-0301122.

\bibliographystyle{amsplain}
\bibliography{tril}

\end{document}